\documentclass[preprint,12pt]{elsarticle}

\usepackage{amssymb}
\usepackage{amsmath}
\usepackage{graphicx}
\usepackage[usenames,dvipsnames,table]{xcolor}
\usepackage{amsfonts}
\usepackage{marginnote}
\usepackage{color}
\usepackage{bm}
\usepackage{cancel}
\usepackage{url}
\usepackage{oldgerm}
\usepackage{wrapfig}
\usepackage{lscape}
\usepackage{rotating}
\usepackage{wasysym}
\usepackage{soul}
\usepackage{amsthm}
\usepackage{verbatim}
\usepackage{subcaption}
\usepackage{tabularx}
\usepackage{booktabs}
\usepackage{siunitx}

\usepackage{tikz}
\usetikzlibrary{quotes,arrows.meta}
\usepackage{pgfplots}

\newtheorem{remark}{Remark}

\newcommand{\mcJ}{\mathcal{J}}

\newcommand{\mcF}{\mathcal{F}}

\newcommand{\mcL}{\mathcal{L}}

\newcommand{\mbR}{\mathbb{R}}
\newcommand{\mbRd}{{\mathbb{R}^d}}

\newcommand{\Omg}{{\Omega}}
\newcommand{\omg}{{\omega}}
\newcommand{\etac}{{\eta_c}}
\newcommand{\etan}{{\eta_n}}
\newcommand{\etad}{{\eta_D}}
\newcommand{\Gamc}{{\Gamma_c}}
\newcommand{\Gaml}{{\Gamma_l}}
\newcommand{\Gamd}{{\Gamma_D}}
\newcommand{\omgn}{{\omega_n}}
\newcommand{\Omgn}{{\Omega_n}}
\newcommand{\Omgl}{{\Omega_l}}
\newcommand{\Omgo}{{\Omega_o}}

\def \bb{\mathbf{b}}
\def \gb{\mathbf{g}}
\def \nub{{\boldsymbol \nu}}
\def \fb{\mathbf{f}}
\def \ub{\mathbf{u}}

\def \ubn{\ub_n}
\def \ubl{\ub_l}

\def \Ubn{\maxvec{\bm{u}}_n}
\def \Ubl{\maxvec{\bm{u}}_l}
\def \Ubb{\maxvec{\bm{u}}}
\def \Un{\maxvec{u}_n}
\def \Ul{\maxvec{u}_l}
\def \nubn{\nub_n}
\def \nubl{\nub_l}
\def \xb{\bm{x}}
\def \xbp{{\bm{x}^\prime}}
\def \yb{\bm{y}}
\def \Tb{\mathbf{T}}

\def \xib{{\boldsymbol\xi}}
\def \zetab{{\boldsymbol\zeta}}

\def \zerob{\mathbf{0}}

\DeclareMathAccent{\maxvec}{\mathord}{letters}{"7E}

\def\CC{{C\nolinebreak[4]\hspace{-.05em}\raisebox{.4ex}{\tiny\bf ++}}}

\journal{XXX}

\begin{document}
\begin{frontmatter}

\title{An optimization-based strategy for peridynamic-FEM coupling and for the prescription of nonlocal boundary conditions\footnote{Sandia National Laboratories is a multimission laboratory managed and operated by National Technology and Engineering Solutions of Sandia, LLC., a wholly owned subsidiary of Honeywell International, Inc., for the U.S. Department of Energy's National Nuclear Security Administration under contract DE-NA-0003525.}}

\author{Marta D'Elia}

\address{Computational Science and Analysis\\
         Sandia National Laboratories \\
         Livermore, CA}
         
\author{David J.~Littlewood}
\author{Jeremy Trageser}
\author{Mauro Perego}
\author{Pavel B.~Bochev}

\address{Center for Computing Research\\
         Sandia National Laboratories \\
         Albuquerque, NM}

\begin{abstract}
We develop and analyze an optimization-based method for the coupling of a static peridynamic (PD) model and a static classical elasticity model.  The approach formulates the coupling as a control problem in which the states are the solutions of the PD and classical equations, the objective is to minimize their mismatch on an overlap of the PD and classical domains, and the controls are virtual volume constraints and boundary conditions applied at the local-nonlocal interface. Our numerical tests performed on three-dimensional geometries illustrate the consistency and accuracy of our method, its numerical convergence, and its applicability to realistic engineering geometries.  We demonstrate the coupling strategy as a means to reduce computational expense by confining the nonlocal model to a subdomain of interest, and as a means to transmit local (e.g., traction) boundary conditions applied at a surface to a nonlocal model in the bulk of the domain.
\end{abstract}


\end{frontmatter}

\section{Introduction} \label{sec:intro}
Nonlocal models have become viable alternatives to classical partial differential equation (PDE) models for certain classes of problems, particularly those in which discontinuities such as cracks are present. The use of these models is increasing in several scientific and engineering applications such as fracture mechanics \cite{Ha2011,Silling2000}, anomalous subsurface transport \cite{Benson2000,Gulian2021,Schumer2003,Schumer2001}, phase transitions \cite{Burkovska2021,Delgoshaie2015,Fife2003}, image processing \cite{Buades2010,DElia2021Imaging,Gilboa2007}, stochastic processes \cite{Burch2014,DElia2017,Meerschaert2012,MeKl00}, and turbulence \cite{DiLeoni2021,Pang2020npinns,Pang2019fPINNs}. In this work we specifically focus on mechanics applications modeled by peridynamics, a nonlocal extension of continuum mechanics developed to capture discontinuities that result from material failure, as well as other phenomena \cite{Silling_00_JMPS,Silling05,Silling_07_JE}. However, several applications mentioned above, especially the ones where long-range interactions are limited to finite regions, may benefit from the technique presented in this work.

\smallskip
Spatial nonlocal operators are integral operators that embed length scales in their definition; the most general form of a nonlocal operator acting on a vector function $\ub:\mbRd\to\mbRd$ is given by
$$
\mcL[\ub](\xb) = \int_{B_\delta(\xb)} \fb(\xb,\yb;\ub)\,d\yb
$$
where $B_\delta(\xb)$ is the ball centered at $\xb$ of radius $\delta$, usually referred to in the peridynamics literature as the {\it horizon} or interaction radius. The horizon determines the extent of the nonlocal interactions and embeds the characteristic length scale of the system. Furthermore, since the integrand does not depend on the derivatives of $\ub$, the regularity requirements on the solutions of nonlocal problems are minimal (as opposed to the case of PDEs).

\smallskip
While the integral form allows one to capture multiscale behavior and discontinuities, it also poses theoretical and numerical challenges. Theoretical challenges include the lack of a complete nonlocal theory \cite{Defterli2015,DElia2020Helmholtz,DElia2020Unified}, the nontrivial treatment of nonlocal interfaces \cite{Alali2015,Capodaglio2020,Seleson2013,fan2021asymptotically} and nonlocal boundary conditions \cite{DEliaNeumann2019,DElia2021prescription,Foss2021,You_2019,yu2021asymptotically}, which must be prescribed in a volumetric region surrounding the domain of interest to guarantee the uniqueness of the solution. Computational challenges are related to the integral form that may require sophisticated quadrature rules, yielding discretized systems whose matrices are dense or even full. Among the works dedicated to improving implementation of nonlocal discretizations and numerical solvers for discretized nonlocal problems we mention variational methods \cite{Aulisa2021,Capodaglio2021DD,DElia-ACTA-2020,DEliaFEM2020} and meshfree methods \cite{Pasetto2019,Silling2005meshfree,trask2019asymptotically,Wang2010,XuFETI2021}.

Despite the recent tremendous effort towards improving the efficiency of nonlocal simulations, the usability of nonlocal models is still hindered by their computational cost. This work strives to improve the viability of nonlocal methods for systems in which nonlocal modeling can be restricted to certain subregions within the overall domain. It is often the case that in engineering mechanics problems, the application of a nonlocal model is advantageous only in specific regions of a body, for example due to the presence of cracks. In these circumstances it is natural to model the material displacement with a nonlocal description in the proximity of the discontinuity and to use classical PDE models elsewhere, so that the bulk of the computation is concentrated only on select nonlocal regions. Furthermore, by using local models far from these regions, many classical codes that contain an array of features not widely available in peridynamics codes (such as structural elements) can be exploited.

The approach of unifying local and nonlocal models within a single computational framework is known in the literature as {\it Local-to-Nonlocal (LtN) coupling} and has been the subject of very active research during the last decade. According to the review paper \cite{DElia2021review}, LtN approaches can be divided in two categories: constant horizon approaches, where the nonlocal region is characterized by a constant value of the horizon and the transition to the local region is abrupt, and variable horizon approaches, characterized by a smooth transition from the nonlocal to the local region by means of a varying horizon. Among the first category we mention Optimization-Based Methods (OBM) \cite{DElia_14_INPROC,DEliaCoupling,delia2016,DElia2021coupling}, partitioned procedures \cite{You2020,yu2018partitioned}, Arlequin approaches \cite{Han_12_IJNME,ArlequinWang2019}, morphing methods \cite{han2016morphing,Lubineau_12_JMPS}, quasi-nonlocal methods \cite{DuLiLuTian2018,XHLiLu2017}, blending methods \cite{Seleson_13_CMS,Seleson15}, and splice methods \cite{galvanetto2016effective}. Among the second category we mention shrinking horizon methods \cite{TTD19} and partial stress methods \cite{silling2015variable}. A common feature of all these methods is the fact that the domain of interest is divided into possibly overlapping nonlocal and local regions under the assumption that there exists a local model that accurately describes the behavior of the system when the nonlocal effects vanish. As we  describe in more detail below, the most widely used nonlocal models for mechanics converge (pointwise or in norm) to a well-known PDE model \cite{Seleson15} as the extent of the nonlocal interaction vanishes, i.e., as $\delta\to 0$.

\smallskip
Building on our previous work on OBM for nonlocal diffusion problems \cite{DElia_14_INPROC,delia2016,DElia2021coupling} and our preliminary results on OBM for mechanics \cite{DEliaCoupling}, in this work, we propose a physically consistent and non-intrusive OB coupling scheme for the coupling of a state-based peridynamic model and a classical, linear elasticity model. Specifically, we address the LtN coupling of the Linear Peridynamic Solid (LPS) model \cite{Silling_07_JE} and the Navier-Cauchy classical model. The main feature of OBMs, that distinguishes these methods from other LtN techniques, is the fact that coupling conditions between local and nonlocal models are prescribed weakly, via optimization, while the local and nonlocal equations act as constraints of the optimization problem. More specifically, OBMs are formulated as optimization problems where the objective functional is the norm of the mismatch between local and nonlocal displacements in the overlap of the local and nonlocal subdomains, the constraints are the local and nonlocal equations in their respective subdomains, and the control variables are the local boundary condition and nonlocal volume constraint on the fictitious interfaces induced by the domain decomposition. 

We summarize the main advantages OBMs below.
\begin{enumerate}
    \item The constraints are fully decoupled, i.e.~they can be solved independently, as opposed to other coupling methods such as blending approaches that introduce hybrid models via linear combination of local and nonlocal forces or energies. As a consequence, OBMs may be applied either within a unified computational framework, or within a framework in which the local and nonlocal models are evaluated independently, as so-called black boxes.\label{adv:decoupling}
    \item An immediate consequence of \ref{adv:decoupling} is that nonlocal and local discretizations, grids, and software can be different and independent. In fact, the only firm requirement is that the nonlocal and local solutions be compared on the overlap region, such that a simple evaluation or projection operator suffices.  Additional information, for example evaluation of gradients, are optional depending on the nature of the optimization strategy employed. These facts make OBMs extremely flexible in terms of implementation. 
    \item OB solutions inherit the mathematical and numerical properties of the constraints. More specifically: when the constraints are well-posed problems, the OB formulation is also well-posed. Moreover, the overall numerical convergence properties of the OB solutions are the same as those of the schemes used for the discretization of the constraints (as illustrated by three-dimensional numerical tests in the current study). 
    \item OBMs for LtN coupling are amenable to the conversion of local (e.g., surface) boundary conditions to nonlocal volume constraints. It is often the case that volume constraints, necessary for the well-posedness of nonlocal equations, are not readily available. By placing a local model in the vicinity of the boundary, one can use available surface data as boundary conditions for the local model and utilize the coupling scheme to transmit the effect of the boundary conditions to the nonlocal domain. A demonstration of this technique is presented in Section \ref{numerical-tests}.
\end{enumerate}

The major contribution of this manuscript is consolidation of OB coupling methods for mechanics problems in three dimensions by demonstrating their consistency and convergence properties and by illustrating their applicability on representative test cases. Additionally, we show how to employ this technique to facilitate the prescription of boundary conditions in the absence of volumetric (nonlocal) data.

\paragraph{Outline of the paper} In Section \ref{sec:state-models} we describe the peridynamic model and its classical counterpart, i.e.~the classical linear elasticity equation. In Section \ref{sec:optimization} we introduce the OB coupling scheme and highlight some of its properties. In Section \ref{sec:discretization} we provide details on the discretization techniques used in this work, namely, a meshfree approach for the LPS model and the finite element method for the classical model. Finally, in Section \ref{numerical-tests} we illustrate the consistency, accuracy, and applicability of our scheme using numerical examples.

\section{The peridynamics model and its local counterpart} \label{sec:state-models}

In this section we introduce the static, linearized peridynamic equilibrium equation \cite{Silling_07_JE} and its corresponding local limit, i.e.~the classical Navier-Cauchy equation of static elasticity. 

Given the bounded body $\Omg\subset\mbR^3$ with boundary $\partial\Omg=\Gamma$, the general peridynamic equation of the motion at the material point $\xb\in\Omg$ at time $t\geq 0$ is given by
\begin{displaymath}
\rho(\xb,t)\;\frac{\partial^2\ub}{\partial t^2}(\xb,t) = \!\int_\Omg \left\{ \Tb[\xb,t]\langle\xbp-\xb\rangle - \Tb[\xbp,t]\langle\xb-\xbp\rangle \right\} d V_\xbp + \bb(\xb,t).
\end{displaymath}
Here, 
$\rho\!:\Omg\!\times\!\mbR^+\!\!\rightarrow\!\mbR^+$ is the mass density, 
$\ub\!:\Omg\!\times\!\mbR^+\!\!\rightarrow\!\mbR^3$ is the displacement field, 
$\bb\!:\Omg\!\times\!\mbR^+\!\!\rightarrow\!\mbR^3$ is a given body force density and 
$\Tb\!:\Omg\!\times\!\mbR^+\!\!\rightarrow\!\mbR^{(3,3)}$ is the force state field. The latter represents the force state at $(\xb,t)$ that maps the bond $\langle\xbp-\xb\rangle$ to force per unit volume squared. For simplicity, in this work we consider the static problem, for which the equilibrium equation reads
\begin{equation}\label{eq:PDstatic}
-\mcL_n[\ub](\xb):=- \!\int_\Omg \left\{ \Tb[\xb]\langle\xbp-\xb\rangle - 
\Tb[\xbp]\langle\xb-\xbp\rangle \right\} d V_\xbp = \bb(\xb).
\end{equation}
A fundamental assumption in nonlocal modeling is that a material point $\xb$ interacts with material points in a ball centered in $\xb$ of radius $\delta$, which we refer to as the {\it horizon} or nonlocal interaction radius. Formally, this nonlocal neighborhood is given by 
\begin{displaymath}
B_\delta(\xb) = \{\xbp\in\Omg: |\xb-\xbp|\leq\delta\},
\end{displaymath}
As a consequence, we require that the force state field satisfies the following property:
\begin{displaymath}
\Tb[\xb]\langle\xbp-\xb\rangle = 0, \;\; \forall\,\xbp \notin B_\delta(\xb).
\end{displaymath}

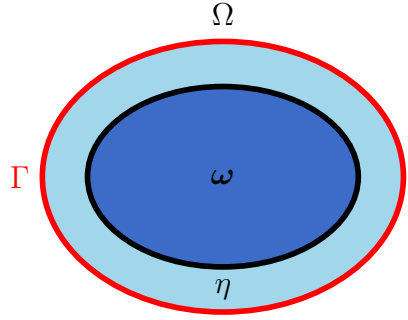
\begin{figure}
\begin{center}
\begin{tikzpicture}[scale=0.6]
\draw[fill={rgb,255:red,162; green,214; blue,235}, draw=red, thick,line width=0.75mm] (0,0) ellipse (4.0cm and 3.0cm);

\draw[fill={rgb,255:red,60; green,108; blue,199}, draw=black, thick,line width=0.75mm] (0,0) ellipse (3.0cm and 2.0cm);

\node[black] at (0.0,0.0) { $\boldsymbol{\omega}$};
\node[black] at (0.0cm,-2.5cm) { $\eta$};
\node[red] at (-4.5cm,0.0cm) { $\Gamma$};
\node[black] at (0.0cm,3.6cm) { $\Omega$};
\end{tikzpicture}
\caption{A domain configuration in two-dimensions.}\label{fig:full-domain}
\end{center}
\end{figure}

Equation \eqref{eq:PDstatic} is solved in the body $\omg\!\subset\!\Omg$ and Dirichlet volume constraints are prescribed in a volumetric layer $\eta$ surrounding $\omg$. These domains are such that $\Omg=\omg\cup\eta$, see Figure \ref{fig:full-domain} for a two-dimensional illustration. The thickness of $\eta$ depends on the definition of $\Tb$ and guarantees, together with other conditions, that the peridynamics problem is well-posed. According the Linear Peridynamic Solid (LPS) model considered in this work \cite{Silling_07_JE}, for all $\xb\in\Omg$, the force state field $\Tb$ is defined as
\begin{equation}\label{eq:linLPS}
\Tb[\xb]\langle\xib\rangle = \dfrac{\kappa(|\xib|)}{m}\left\{ 
(3K-5G)\,\theta(\xb)\xib+15G\frac{\xib\otimes\xib}{|\xib|^2}(\ub(\xb+\xib)-\ub(\xb))\right\},
\end{equation}
where $\xib=\xbp-\xb$, $K$ is the bulk modulus, $G$ is the shear modulus, and the linearized nonlocal dilatation, $\theta\!:\!\Omg\!\rightarrow\!\mbR$, is defined as
\begin{displaymath}
\theta(\xb)= \frac{3}{m}\int_{B_\delta({\bf 0})}\kappa(|\zetab|)\,\zetab\!\cdot\!(\ub(\xb+\zetab)-\ub(\xb))\,dV_{\zetab},
\;\;\hbox{with}\;\;
m = \int_{B_\delta({\bf 0})}\kappa(|\zetab|)\,|\zetab|^2\,dV_{\zetab}.
\end{displaymath}
The spherical function $\kappa$ is referred to as the {\it influence function}; it determines the support of force states and modulates the bond strength \cite{Seleson_11,Silling_07_JE}. With this choice of $\Tb$, we define the volumetric layer $\eta$ as
\begin{equation}\label{eq:layer-definition}
\eta=\{\xbp\in\Omg: \; |\xbp-\xb|< 2\delta\} \quad \forall \, \xb\in \Gamma,
\end{equation}
whose thickness is double the size of the horizon. As mentioned before, this choice guarantees that the peridynamic operator can be evaluated for any point in $\omg$, up to its boundary $\partial\omg$. This is due to the fact that $\theta$ introduces an additional integral in the definition of $\mcL_n$; thus, we have a double integral over $B_\delta(\zerob)\!\times\! B_\delta(\zerob)$, i.e.~for every point $\xb$, we need values of the displacement in $B_{2\delta}(\xb)$.

Using the linearized LPS force state field in \eqref{eq:linLPS} the three-dimensional peridynamic equation is then formulated as follows. Find $\ub\in [L^2(\Omg)]^3$ such that
\begin{equation}\label{eq:linLPS-problem}
\left\{
\begin{array}{rcl}
- \mcL_{\rm LPS}[\ub](\xb) =& \!\bb(\xb) & \;\; \xb\in\omg \\[3mm]
\ub(\xb) = & \!\gb(\xb) &\;\; \xb\in\eta.
\end{array}\right.
\end{equation}
Here, $\gb\in[L^2(\eta)]^3$ is the volumetric Dirichlet data and $\mcL_{\rm LPS}$ is the nonlocal operator $\mcL_n$ corresponding to the choice of $\Tb$ as in \eqref{eq:linLPS}, i.e.
\begin{equation}\label{eq:LLPS}
\begin{aligned}
\mcL_{\rm LPS}[\ub](\xb)
& := \int_{B_\delta({\bf 0})} \dfrac{\kappa(|\xib|)}{m}
\left\{\phantom{\dfrac{1}{1}}\!\!\!\!  (3K-5G)(\theta(\xb)+\theta(\xb+\xib))\xib\right. \\[3mm]
& \; + 30G\left.\frac{\xib\otimes\xib}{|\xib|^2}(\ub(\xb+\xib)-\ub(\xb))\right\}\,dV_{\xib}.
\end{aligned}
\end{equation}

The linearized LPS model \eqref{eq:LLPS} has two important properties that one can exploit to choose a viable candidate for the local model and perform consistency and converge tests for the coupling method. First, the limit for $\delta\to 0$ (i.e., when the nonlocal interactions vanish) is the classical Navier-Cauchy equation (NCE) of static elasticity \cite{Seleson2016}:
\begin{equation}\label{eq:NCE}
- \mcL_{\rm NC}[\ub](\xb) := -\left[\left(K+\frac13 G\right)\nabla(\nabla\cdot \ub)(\xb)+G\,\nabla^2\ub(\xb)\right]= \bb(\xb),
\end{equation}
where $K$, $G$ and $\bb$ are defined as in \eqref{eq:linLPS}, or equivalently,
\begin{equation}\label{eq:NCElame}
\begin{array}{l}
\displaystyle-\nabla\cdot{\boldsymbol\sigma}[\ub](\xb)  =  \bb(\xb), \quad {\rm where}\\[4mm]
{\boldsymbol\sigma}[\ub] =  \lambda (\nabla\cdot\ub) \mathbf I +
                            \mu(\nabla\ub + \nabla\ub^T)\\[3mm]
(\lambda,\mu) =  \left(K - \frac{2G}{3}, \, G\right),
\end{array}
\end{equation}
where $\mathbf I$ is the identity tensor. Note that the latter equation is the classical linear elasticity equation in terms of the Lam\'e constants $(\lambda,\mu)$. This implies that the NC model is a fair approximation of the peridynamics model for sufficiently regular solutions and, as a consequence, it can be used as the local model in the coupling strategy. 

Furthermore, the linearized LPS model equates to the classical NCE model when the displacement field is quadratic~\cite[Proposition 1]{Seleson2016}. This property provides analytical solutions for the LPS model for the quadratic patch test and numerical convergence tests presented in Section~\ref{sec:patch-convergence}.

\section{Optimization-based LtN formulation}\label{sec:optimization}
\begin{figure}
\begin{center}
\begin{tikzpicture}[scale=0.8]

\def\OmegaLColor{rgb,255:red,60; green,108; blue,199};
\def\omeganColor{rgb,255:red,162; green,214; blue,235};
\def\etacColor{yellow};
\def\etaDColor{pink};
\def\GammacColor{black!10!red};
\def\GammaDColor{brown};
\def\horzl{(0,0) ellipse (5.0cm and 3.0cm)};
\def\horzs{(0,0) ellipse (4.5cm and 2.5cm)};
\def\vertl{(0,0) ellipse (1.75cm and 5cm)};
\def\verts{(0,0) ellipse (1.25cm and 4.5cm)};

\def\boxlenhor{6};
\def\boxlenver{3};
\def\leftplane{(-\boxlenhor,-\boxlenver) -- (-\boxlenhor,\boxlenver) -- (0,\boxlenver) -- (0,-\boxlenver) -- cycle};
\def\rightplane{(0,-\boxlenver) -- (0,\boxlenver) -- (\boxlenhor,\boxlenver) -- (\boxlenhor,-\boxlenver) -- cycle};
\def\rightplaneshift{(1,-\boxlenver) -- (1,\boxlenver) -- (\boxlenhor,\boxlenver) -- (\boxlenhor,-\boxlenver) -- cycle};

\clip (-8.0,-3.0) rectangle (8.0,3.0); 

\begin{scope}
   \clip \horzs;
   \clip \leftplane;
   \fill[fill={\omeganColor}, even odd rule] \horzs \vertl;
\end{scope}

\begin{scope}
   \clip \horzl;
   \clip \rightplane;
   \fill[fill={\OmegaLColor}, even odd rule] \horzl \vertl;
\end{scope}

\begin{scope}
   \clip \leftplane;
   \clip \horzl;
   \fill[fill={\etaDColor}, even odd rule] \horzl \horzs;
\end{scope}

\begin{scope}
   \clip \vertl;
   \clip \horzl;

   \foreach \x in {-10,...,20}
    {\pgfmathsetmacro{\a}{-5+2*\x/4}
     \pgfmathsetmacro{\b}{-5/4-5/4*\a}
     \pgfmathsetmacro{\c}{-5+2*\x/4+1/4}
     \pgfmathsetmacro{\d}{-5/4-5/4*\c}
     \fill[fill={\OmegaLColor}] (\a,-5) -- (3,\b) -- (3,-5) -- cycle;
     \fill[fill={\etaDColor}] (\c,-5) -- (3,\d) -- (3,-5) -- cycle;}
\end{scope}

\begin{scope} 
   \clip \horzs;
   \clip \vertl;

   \foreach \x in {-10,...,20}
    {\pgfmathsetmacro{\a}{-5+2*\x/4}
     \pgfmathsetmacro{\b}{-5/4-5/4*\a}
     \pgfmathsetmacro{\c}{-5+2*\x/4+1/4}
     \pgfmathsetmacro{\d}{-5/4-5/4*\c}
     \fill[fill={\OmegaLColor}] (\a,-5) -- (3,\b) -- (3,-5) -- cycle;
     \fill[fill={\omeganColor}] (\c,-5) -- (3,\d) -- (3,-5) -- cycle;}

\end{scope}

\begin{scope} 
   \clip \horzs;
   \clip \vertl;
   \clip \rightplane \verts; 
   \clip \rightplaneshift; 

   \foreach \x in {-10,...,20}
    {\pgfmathsetmacro{\a}{-5+2*\x/4}
     \pgfmathsetmacro{\b}{-5/4-5/4*\a}
     \pgfmathsetmacro{\c}{-5+2*\x/4+1/4}
     \pgfmathsetmacro{\d}{-5/4-5/4*\c}
     \fill[fill={\OmegaLColor}, even odd rule] {(\a,-5) -- (3,\b) -- (3,-5) -- cycle} {\verts};
     \fill[fill={\etacColor}, even odd rule] {(\c,-5) -- (3,\d) -- (3,-5) -- cycle} {\verts};}

\end{scope}


\begin{scope}
   \clip \horzs;
   \clip \vertl;
   \draw[black, line width=1.0mm] \horzs;
\end{scope}
\begin{scope}
   \clip \horzs;
   \clip \leftplane;
   \draw[black, line width=1.0mm] \horzs;
\end{scope}
\begin{scope}
   \clip \horzl;
   \clip \rightplane;
   \draw[black, line width=0.5mm] \vertl;
\end{scope}
\begin{scope}
   \clip \horzs;
   \clip \rightplane;
   \draw[black, line width=0.5mm] \verts;
\end{scope}
\begin{scope}
   \clip \vertl;
   \clip \horzl;
   \clip \leftplane;
   \draw[color={\GammacColor}, line width=1.25mm] \vertl;
\end{scope}
\begin{scope} 
   \clip \horzl;
   \clip (-1.43193,-10) -- (-1.43193,10) -- (10,10) -- (10,-10) -- cycle;
   \draw[color={\GammaDColor}, line width=1.0mm] \horzl;
\end{scope}
\begin{scope}
   \clip (-1.43193,-10) -- (-1.43193,10) -- (-10,10) -- (-10,-10) -- cycle;
   \clip \horzl;
   \draw[black,line width=1.0mm] \horzl;
\end{scope}


\draw[color={\GammaDColor}, line width=0.5mm] (6,1.5) -- (7,1.5) node[right,text=black] {$\Gamma_D$};
\draw[color={\GammacColor}, line width=0.5mm] (6,0.25) -- (7,0.25) node[right,text=black] {$\Gamma_c$};
\draw[color={\OmegaLColor},line width=2.5mm] (6,-1.0) -- (7,-1.0) node[right,text=black] {$\Omega_l$};

\draw[color={\omeganColor},line width=2.5mm] (-6,-1.0) -- (-7,-1.0) node[left,text=black] {$\omega_n$};
\draw[color={\etaDColor},line width=2.5mm] (-6,0.25) -- (-7,0.25) node[left,text=black] {$\eta_D$};
\draw[color={\etacColor},line width=2.5mm] (-6,1.5) -- (-7,1.5) node[left,text=black] {$\eta_c$};

\end{tikzpicture}
\caption{An example LtN domain configuration in two-dimensions.}\label{fig:2D-config}
\end{center}
\end{figure}
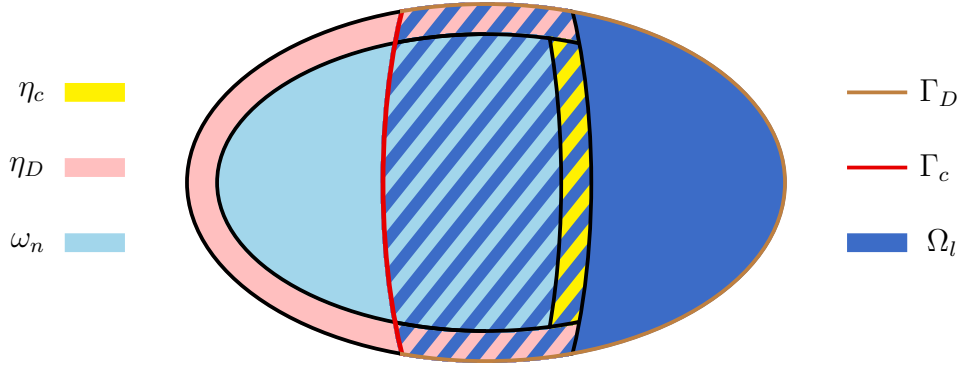

To construct the continuous formulation of the coupling procedure, we first introduce a partitioning of the domain $\Omg$ into a nonlocal subdomain $\Omgn$ and a local subdomain $\Omgl$, with boundary $\Gaml$, such that $\Omgn=\omgn\cup\etan$ and $\Omgn\cap\Omgl=\Omgo\neq\emptyset$. A two-dimensional illustration of this strategy is given in Figure~\ref{fig:2D-config}.
As is common in LtN coupling methods, we assume that the linearized LPS model \eqref{eq:LLPS} provides an accurate description of the material behavior in $\Omgn$ whereas the local NC model gives a fair representation in the remainder of the domain. In OBM approaches to coupling, the main idea is to formulate the coupling as an optimization problem where the difference between the nonlocal and the local displacements is minimized on the overlap $\Omgo$ by tuning their values on the virtual interfaces induced by the partition, i.e.~the virtual interaction volume $\etac$ and the virtual boundary $\Gamc$. Let $\etad=\eta\cap\etan$ and $\Gamd=\Gamma\cap\Gaml$ be the physical interaction volume and boundary where we prescribe the given Dirichlet data. We define the virtual control volume and boundary as $\etac=\etan\!\setminus\!\etad$ and $\Gamc=\Gaml\!\setminus\!\Gamd$. We then pose the peridynamics problem on the nonlocal domain $\omgn$ and the classical NC problem on the local domain $\Omgl$. This results in the following systems of equations.
\begin{equation}\label{eq:constraints}
\left\{
\begin{array}{rcl}
- \mcL_{\rm LPS}[\ubn](\xb) =& \!\bb(\xb) & \xb\in\omgn \\[3mm]
\ubn(\xb) = & \!{\bf g}(\xb) & \xb\in\etad \\[3mm]
\ubn(\xb) = & \!\nubn(\xb) & \xb\in\etac
\end{array}\right.
\quad
\left\{\begin{array}{rcl}
- \mcL_{\rm NC}[\ubl](\xb) =& \!\bb(\xb) & \xb\in\Omgl \\[3mm]
\ubl(\xb) = & \!{\bf g}(\xb) & \xb\in\Gamd \\[3mm]
\ubl(\xb) = & \!\nubl(\xb) & \xb\in\Gamc.
\end{array}\right.
\end{equation}
Here, $\nubn\in[L^2(\etac)]^3$ and $\nubl\in[H^{1/2}(\Gamc)]^3$ are the unknown volume constraint and boundary condition. The systems \eqref{eq:constraints} act as constraints to the optimization problem, while $(\nubn,\nubl)$ act as control variables. As a result, we aim at solving the following optimization problem:
\begin{equation}\label{eq:continuous-opt}
\min\limits_{\ubn,\ubl,\nubn,\nubl} \mcJ(\ubn,\ubl)=\frac12 \int_\Omgo 
|\ubn-\ubl|^2 \,d\xb \quad \hbox{subject to \eqref{eq:constraints}.} 
\end{equation}
Upon solution of \eqref{eq:continuous-opt}, we denote the optimal controls by $\nubn^*$ and $\nubl^*$, and the corresponding optimal displacements by $\ubn^*=\ubn(\nubn^*)$ and $\ubl^*=\ubl(\nubl^*)$. Then, the optimal coupled solution is defined as
\begin{equation}\label{eq:coupled_sol}
\ub^* = \left\{
\begin{array}{ll}
\ubn^* & \xb\in\Omgn \\[2mm]
\ubl^* & \xb\in\Omgl\setminus\Omgo.
\end{array}\right.
\end{equation}

\begin{remark}
The well-posedness of the minimization  \eqref{eq:continuous-opt} has been rigorously studied in \cite{DElia2021coupling} and \cite{delia2016} for nonlocal Poisson's problems with Dirichlet and Neumann conditions in a multi-dimensional setting. Further considerations on the well-posedness in the context of peridynamics models can be found in \cite{DEliaCoupling}.  
\end{remark}

\section{Numerical solution of the optimization-based LtN formulation} \label{sec:discretization}
For the linearized LPS model described in Section \ref{sec:state-models} we utilize the meshfree approach of Silling and Askari~\cite{Silling05}. We discretize the nonlocal domain using a set of material points $\{\xb_i\}_{i=1}^{N_n}\subset\Omg_{n}$, and associate with each material point $\xb_i$ a volume $V_i$ such that $\sum_{i=1}^{N_n} V_i = \Omg_{n}$.  For every point $\xb_i$ we approximate the operator $\mcL_{\rm LPS}$ as follows
\begin{equation}\label{eq:approx-operator}
\mcL_{\rm LPS}^h[\xb_i]:= \sum_{j\in\mcF_i}\left\{\Tb[\xb_i]\langle\xb_j-\xb_i\rangle - 
\Tb[\xb_j]\langle\xb_i-\xb_j\rangle\right\} \, V_j^{(i)},
\end{equation}
where $\xb_i$ and $V_j^{(i)}$ serve as quadrature points and weights and $\mcF_i$ represents the {\it family} of $\xb_i$, i.e., the set of all points $\xb_j$ in $\Omg_{n}$ that are within a distance of $\delta$ from $\xb_i$. Here, the $\xb_j$ is chosen to coincide with the reference position of the $j$th node and the quadrature weight\footnote{Details regarding the computation of $V_j^{(i)}$ can be found in \cite{Seleson2016}.} $V_j^{(i)}$ is the volume associated $\xb_j$. The vector of degrees of freedom of the discrete nonlocal solution at the material points is denoted by $\Ubn= [\Un^{\,1}, \Un^{\,2}, \Un^{\,3}]$, with $\Un^{\,k}\in\mbR^{N_n}$.

We discretize the NC model in \eqref{eq:NCE} by the finite element (FE) method using continuous piecewise linear basis functions. As it is standard, the FE implementation is not discussed in detail. We denote the vector of degrees of freedom of the local discrete solution by $\Ubl = [\Ul^{\,1}, \Ul^{\,2}, \Ul^{\,3}]$, with $\Ul^{\,k} \in\mbR^{N_l}$ where $N_l$ is the number of degrees of freedom of each spatial component over the FE computational mesh. 

Although the peridynamics equation and the NC equation do not directly interact, the difference of the corresponding solutions must be computed on the overlap region $\Omega_o$. To this end, we introduce the nonlocal and local selection matrices that allow us to compute an approximation of the cost functional. Let $S_n\in\mbR^{N_o,N_n}$ be the matrix that selects the components of $\Un^{\,k}$ in $\Omgo$ and $S_l\in\mbR^{N_o,N_l}$ be the operator that evaluates $\Ul^{\,k}$ at the material points $\{\xb_i\}\in\Omgo$. Formally, given the FE basis functions $\{\phi_j\}_{i=1}^{N_l}$, the selection matrices are defined as
\begin{displaymath}
(S_n)_{ij} := \delta_{ij} \; 
\mbox{ and } \; (S_l)_{ij} := \phi_j(\xb_i), \quad \forall \, \xb_i\in\Omgo.
\end{displaymath}
Then, we define the discrete functional as

\begin{equation}\label{eq:discrete-functional}
J_d(\Ubn,\Ubl) = \frac{1}{2} \sum_{i=1}^{N_o} \sum_{k=1}^{3}|(S_n\Un^{\,k})_i-(S_l\Ul^{\, k})_i|^2 \, \widetilde V_i,
\end{equation}
where $\widetilde V_i$ is the volume associated with the $i$th material point, properly scaled.

\section{Numerical tests}\label{numerical-tests}
In this section we report the results of several numerical tests performed on three-dimensional geometries. The purpose of these experiments is, first, to illustrate the consistency and accuracy of our method via patch tests and numerical convergence studies, and, second, to demonstrate the applicability of our approach for realistic engineering geometries. We also show that our coupling approach can be used to circumvent the nontrivial prescription of ``nonlocal boundary conditions'', or volume constraints, when only surface data are available. We start by providing some details on the software used for the simulations. These are followed by consistency and convergence tests with manufactured solutions. We then consider examples of coupling in the presence of cracks and provide an example of prescribing boundary conditions via OB coupling.

The simulations presented in this work were performed using \emph{Peridigm}~\cite{sand2012-7800,peridigm-website} and \emph{Albany}~\cite{sand2013-8430J,albany-website}, two open-source codes developed primarily at Sandia National Laboratories. \emph{Peridigm} is a peridynamics code for solid mechanics based on the meshfree approach of Silling and Askari \cite{Silling05,sand2015-littlewood}, and \emph{Albany} is a FE code for the simulation of several physical systems governed by PDEs.  Both \emph{Peridigm} and \emph{Albany} are \CC~codes designed for use on large-scale parallel computing platforms. For the current study, \emph{Peridigm} and \emph{Albany} were coupled directly, resulting in a single executable in which the peridynamic and FE contributions are treated in a unified fashion and the implicit solver acts on a single, monolithic system.  Routines in the \emph{Peridigm} code base are used to evaluate peridynamic forces and the corresponding entries in the tangent stiffness matrix. The \emph{Albany} code handles most remaining aspects of the computation, including FE assembly for the NC equation, computation of the discrete functional and its derivative, and solution of the state and adjoint systems necessary for the numerical solution of the optimization problem, conducted by means of the LBFGS algorithm. The \emph{Peridigm} and \emph{Albany} codes, and the software infrastructure developed in this study to couple them, rely heavily on several \emph{Trilinos}~\cite{trilinos-website} packages, including \emph{ROL} for solution of the optimization problem, \emph{Epetra} for the management of parallel data structures, \emph{Intrepid} for FE assembly, and \emph{Ifpack} and \emph{AztecOO} for the preconditioning and solution of the linear systems.

\subsection{Patch tests and convergence study}\label{sec:patch-convergence}
We consider analytic solutions of the coupling problem in order to perform patch tests and a convergence study. We utilize linear displacement fields for the patch test and quadratic displacement fields for the convergence study. As discussed above, it can be shown (see \cite{Seleson2016}) that equations \eqref{eq:linLPS-problem} and \eqref{eq:NCE} are equivalent for quadratic displacements, thus the known analytical solution to the local problem also holds for the nonlocal problem under this condition. In our tests we let $\Omg^+=[0,1]^3$ and $\Omg=[2\delta,1-2\delta]^3$, so that the thickness of the external layer where the volume constraints are prescribed is twice the size of the horizon. Specifically, for $\xb=(x,y,z)$, we consider the following data sets.
\begin{itemize}
\item[I:] linear\\[1mm]
          $\ub_\ell(\xb) = (x,0,0)$, $\bb_\ell(\xb) = {\bf 0}$, $\gb_\ell(\xb) = \ub_\ell(\xb)$.
\item[II:] quadratic\\[1mm]
           $\ub_q(\xb) = (x^2,0,0)$,\\[1mm]
           $\bb_q(\xb) = (-304.49,0,0)$,
           obtained by substitution of $\ub$ in \eqref{eq:NCE},\\
           $\gb_q(\xb) = \ub_q(\xb)$,\\[1mm]
           $(\lambda,\mu) = (109.62,73.08)$ (Lam\'e parameters representative of stainless steel).
\end{itemize}

To assess the accuracy of the coupled solution, we consider the $\ell^2$ norm of the difference between the vectors of the computed displacements and the analytic displacements $\Ubb^*$ at the nonlocal and local degrees of freedom, i.e. 
\begin{equation}\label{eq:error}
\begin{aligned}
{\rm error}_n&=\|\Ubn - \Ubb_*\|_{\ell^2} \\[2mm]
{\rm error}_l&=\|\Ubl - \Ubb_*\|_{\ell^2}
\end{aligned}
\end{equation}
where $\Ubb^*$ is either $\ub_\ell$ or $\ub_q$ evaluated at the appropriate degrees of freedom. 

We perform a linear patch test to verify the consistency of the method; to this end, we prescribe $\bb$ and $\gb$ as in case I above\footnote{For more numerical tests regarding linear and quadratic patch tests in three-dimensional settings we refer the interested reader to \cite{DEliaCoupling}.}. Numerical results indicate a perfect match of the nonlocal and local solutions on the overlap and show that the coupled solution is $\epsilon$-machine accurate with respect to $\ub_\ell$. This happens because linear solutions can be approximated exactly by both the peridynamic and FE discretizations. 

\begin{figure}[t]
\centering
\includegraphics[width=\textwidth]{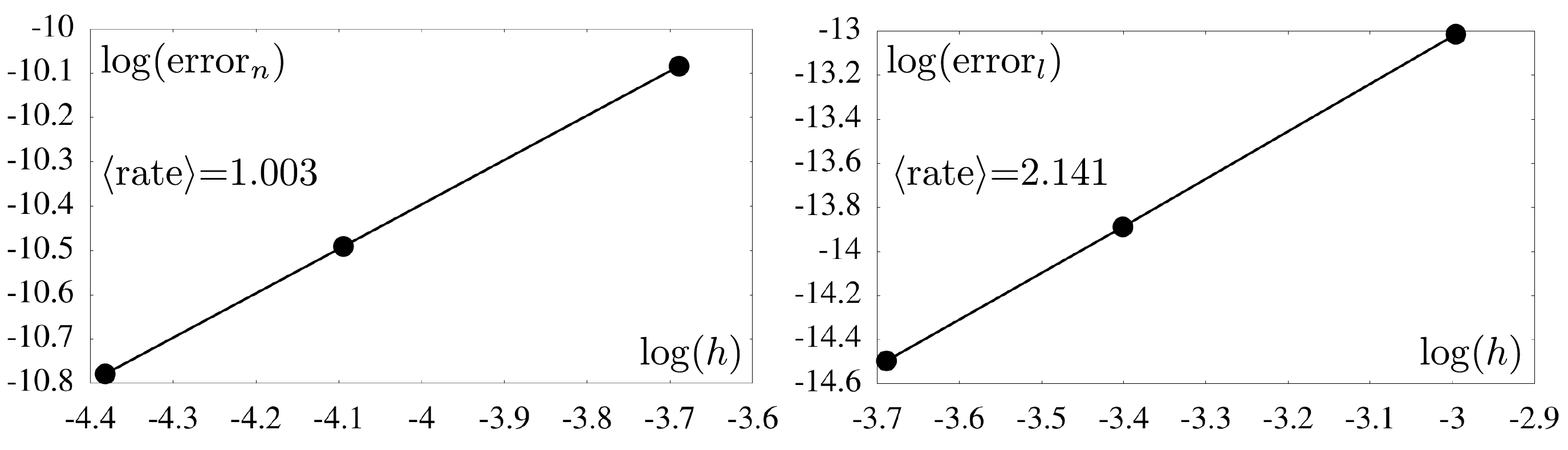}
\caption{Convergence plots for the nonlocal (left) and local (right) portions of the coupled model.}
\label{fig:convergence}
\end{figure}

Results of the convergence tests are reported in Figure \ref{fig:convergence}. Here, we prescribe $\bb$ and $\gb$ as in case II above. For decreasing values of $h$, we report the errors corresponding to the nonlocal and local solutions, according to \eqref{eq:error}, using a log-log scale. We observe a linear rate of convergence for the nonlocal solution and a quadratic rate of convergence for the local solution, see the average rate, $\langle{\rm rate}\rangle$, reported in the figure. This is in line with theoretical results for the meshfree nonlocal discretization utilized in the nonlocal region \cite{Seleson2019volume,Seleson2016}, and for the piecewise linear finite element approximation utilized in the local region \cite{Chen2011,DElia2021coupling}.

%
\subsection{Geometries with cracks}\label{crack}
Peridynamic models are advantageous for modeling cracks because the governing equations remain valid in the presence of material discontinuities.  We investigate the performance of the OB coupling scheme for this class of problems using a three-dimensional bar with a crack. We assign material parameters that are representative of steel alloys, as shown in Table \ref{table::DirichletDirichletBarDomain}.
\begin{table}[ht]
\begin{center}
      \begin{tabular}{lll}   
         \toprule
         \rowcolor{black!15} Parameter & Prenotched bar & Compact Tension\\\midrule
         Horizon & $\SI{1.00}{\milli\meter}$ & $\SI{0.30}{\milli\meter}$ \\\midrule
         Poisson's ratio & $0.3$ & $0.27$ \\\midrule
         Bulk modulus & $\SI{140.0}{\giga\pascal}$ & $\SI{150.0}{\giga\pascal}$\\
         \bottomrule  
      \end{tabular}
         \caption{Material parameters for the bar and compact tension experiments.}\label{table::DirichletDirichletBarDomain}  
\end{center}
\end{table}

\subsubsection{Prenotched rectangular bar}\label{sec:prenotched_bar}

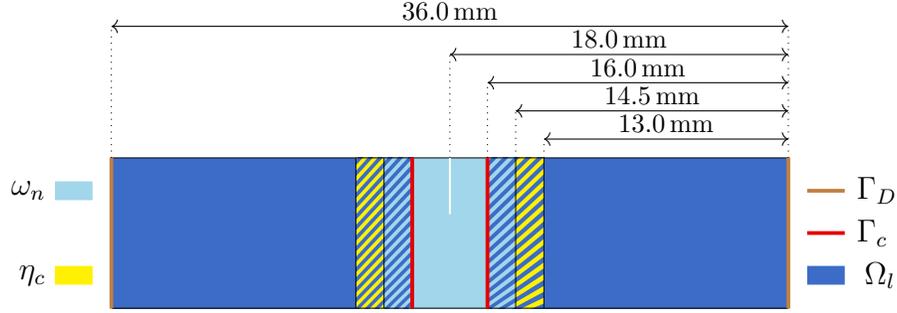
\begin{figure}
\begin{center}
\begin{tikzpicture}[scale=0.25]

\def\OmegaLColor{rgb,255:red,60; green,108; blue,199};
\def\omeganColor{rgb,255:red,162; green,214; blue,235};
\def\etacColor{yellow};
\def\etaDColor{pink};
\def\GammacColor{black!10!red};
\def\GammaDColor{brown};

\def\omeganregion{(-2.0cm,-4.0cm) rectangle (2.0cm,4.0cm)};
\def\omegalleftregion{(-18cm,-4.0cm) rectangle (-5.0cm,4.0cm)};
\def\omegalrightregion{(5.0cm,-4.0cm) rectangle (18.0cm,4.0cm)};
\def\omegalomeganregionl{(-3.5cm,-4.0cm) rectangle (-2.0cm,4.0cm)};
\def\omegalomeganregionr{(2.0cm,-4.0cm) rectangle (3.5cm,4.0cm)};
\def\omegaletacregionl{(-5.0cm,-4.0cm) rectangle (-3.5cm,4.0cm)};
\def\omegaletacregionr{(3.5cm,-4.0cm) rectangle (5.0cm,4.0cm)};

\def\horzs{(0,0) ellipse (4.5cm and 2.5cm)};
\def\vertl{(0,0) ellipse (1.75cm and 5cm)};
\def\verts{(0,0) ellipse (1.25cm and 4.5cm)};

\fill [fill=\omeganColor] \omeganregion;

\fill [fill=\OmegaLColor] \omegalleftregion;
\fill [fill=\OmegaLColor] \omegalrightregion;


\begin{scope} 
  \clip \omegalomeganregionl;

  \foreach \x in {-20,...,20}
    {\pgfmathsetmacro{\a}{-5+2*\x/4}
     \pgfmathsetmacro{\b}{-5/4-5/4*\a}
     \pgfmathsetmacro{\c}{-5+2*\x/4+1/4}
     \pgfmathsetmacro{\d}{-5/4-5/4*\c}
     \fill[fill={\OmegaLColor}] (\a,-5) -- (5,\b) -- (5,-5) -- cycle;
     \fill[fill={\omeganColor}] (\c,-5) -- (5,\d) -- (5,-5) -- cycle;}
\end{scope}
\begin{scope} 
  \clip \omegalomeganregionr;

  \foreach \x in {-20,...,30}
    {\pgfmathsetmacro{\a}{-5+2*\x/4}
     \pgfmathsetmacro{\b}{-5/4-5/4*\a}
     \pgfmathsetmacro{\c}{-5+2*\x/4+1/4}
     \pgfmathsetmacro{\d}{-5/4-5/4*\c}
     \fill[fill={\OmegaLColor}] (\a,-5) -- (5,\b) -- (5,-5) -- cycle;
     \fill[fill={\omeganColor}] (\c,-5) -- (5,\d) -- (5,-5) -- cycle;}
\end{scope}

\begin{scope} 
  \clip \omegaletacregionl;

  \foreach \x in {-20,...,20}
    {\pgfmathsetmacro{\a}{-5+2*\x/4}
     \pgfmathsetmacro{\b}{-5/4-5/4*\a}
     \pgfmathsetmacro{\c}{-5+2*\x/4+1/4}
     \pgfmathsetmacro{\d}{-5/4-5/4*\c}
     \fill[fill={\OmegaLColor}] (\a,-5) -- (5,\b) -- (5,-5) -- cycle;
     \fill[fill={\etacColor}] (\c,-5) -- (5,\d) -- (5,-5) -- cycle;}
\end{scope}
\begin{scope} 
  \clip \omegaletacregionr;

  \foreach \x in {-20,...,30}
    {\pgfmathsetmacro{\a}{-5+2*\x/4}
     \pgfmathsetmacro{\b}{-5/4-5/4*\a}
     \pgfmathsetmacro{\c}{-5+2*\x/4+1/4}
     \pgfmathsetmacro{\d}{-5/4-5/4*\c}
     \fill[fill={\OmegaLColor}] (\a,-5) -- (5,\b) -- (5,-5) -- cycle;
     \fill[fill={\etacColor}] (\c,-5) -- (5,\d) -- (5,-5) -- cycle;}
\end{scope}

\draw[color={black}, line width=0.1mm] (-5.0,-4.0cm) -- (-5.0,4.0cm);
\draw[color={black}, line width=0.1mm] (5.0,-4.0cm) -- (5.0,4.0cm);
\draw[color={black}, line width=0.1mm] (-3.5,-4.0cm) -- (-3.5,4.0cm);
\draw[color={black}, line width=0.1mm] (3.5,-4.0cm) -- (3.5,4.0cm);
\draw[color={black}, line width=0.1mm] (-18.0,-4.0cm) -- (18.0,-4.0cm);
\draw[color={black}, line width=0.1mm] (-18.0,4.0cm) -- (18.0,4.0cm);

\draw[color={white}, line width=0.25mm] (0.0,1.0cm) -- (0.0,4.02cm);

\draw[color={\GammaDColor}, line width=0.5mm] (-18.0,-4.02cm) -- (-18.0,4.02cm);
\draw[color={\GammaDColor}, line width=0.5mm] (18.0,-4.02cm) -- (18.0,4.02cm);

\draw[color={\GammacColor}, line width=0.5mm] (-2.0,-4.02cm) -- (-2.0,4.02cm);
\draw[color={\GammacColor}, line width=0.5mm] (2.0,-4.02cm) -- (2.0,4.02cm);

 \draw[color={\GammaDColor}, line width=0.5mm] (19.0,2.25) -- (21.0,2.25) node[right,text=black] {$\Gamma_D$};
 \draw[color={\GammacColor}, line width=0.5mm] (19.0,0.0) -- (21.0,0.0) node[right,text=black] {$\Gamma_c$};
 \draw[color={\OmegaLColor},line width=2.5mm] (19.0,-2.25) -- (21.0,-2.25) node[right,text=black] {$\Omega_l$};

 \draw[color={\omeganColor},line width=2.5mm] (-21.0,2.25) -- (-19.0,2.25) node[left,text=black,xshift=-2ex] {$\omega_n$};
 \draw[color={\etacColor},line width=2.5mm] (-21.0,-2.25) -- (-19.0,-2.25) node[left,text=black,xshift=-2ex] {$\eta_c$};

\draw[<->] (5.0,5.0) -- (18.0,5.0) node[midway, above, rotate=0,yshift=-0.3ex,xshift=0.0ex]{\footnotesize \SI{13.0}{\milli\meter}};
\draw[dotted] (5.0,4.0) -- (5.0,5.0);

\draw[<->] (3.5,6.5) -- (18.0,6.5) node[midway, above, rotate=0,yshift=-0.3ex,xshift=0.0ex]{\footnotesize \SI{14.5}{\milli\meter}};
\draw[dotted] (3.5,4.0) -- (3.5,6.5);

\draw[<->] (2.0,8.0) -- (18.0,8.0) node[midway, above, rotate=0,yshift=-0.3ex,xshift=0.0ex]{\footnotesize \SI{16.0}{\milli\meter}};
\draw[dotted] (2.0,4.0) -- (2.0,8.0);

\draw[<->] (0,9.5) -- (18.0,9.5) node[midway, above, rotate=0,yshift=-0.3ex,xshift=0.0ex]{\footnotesize \SI{18.0}{\milli\meter}};
\draw[dotted] (0.0,4.0) -- (0.0,9.5);

\draw[<->] (-18.0,11.0) -- (18.0,11.0) node[midway, above, rotate=0,yshift=-0.3ex,xshift=0.0ex]{\footnotesize \SI{36.0}{\milli\meter}};
\draw[dotted] (-18.0,4.0) -- (-18.0,11.0);
\draw[dotted] (18.0,4.0) -- (18.0,11.0);

\end{tikzpicture}
\caption{$XY$ cross section of the geometry for the prenotched bar experiment in which Dirichlet boundary conditions are imposed.  The nonlocal model is applied in the vicinity of the prenotch, and the local model is applied elsewhere.}\label{fig:DirichletBarIllustration}
\end{center}
\end{figure}

\begin{figure}
    \begin{center}
    \begin{subfigure}{\textwidth}
     {   \centering
     \includegraphics[scale=0.35,trim=3.1cm 8.1cm 2.1cm 10.3cm,clip]{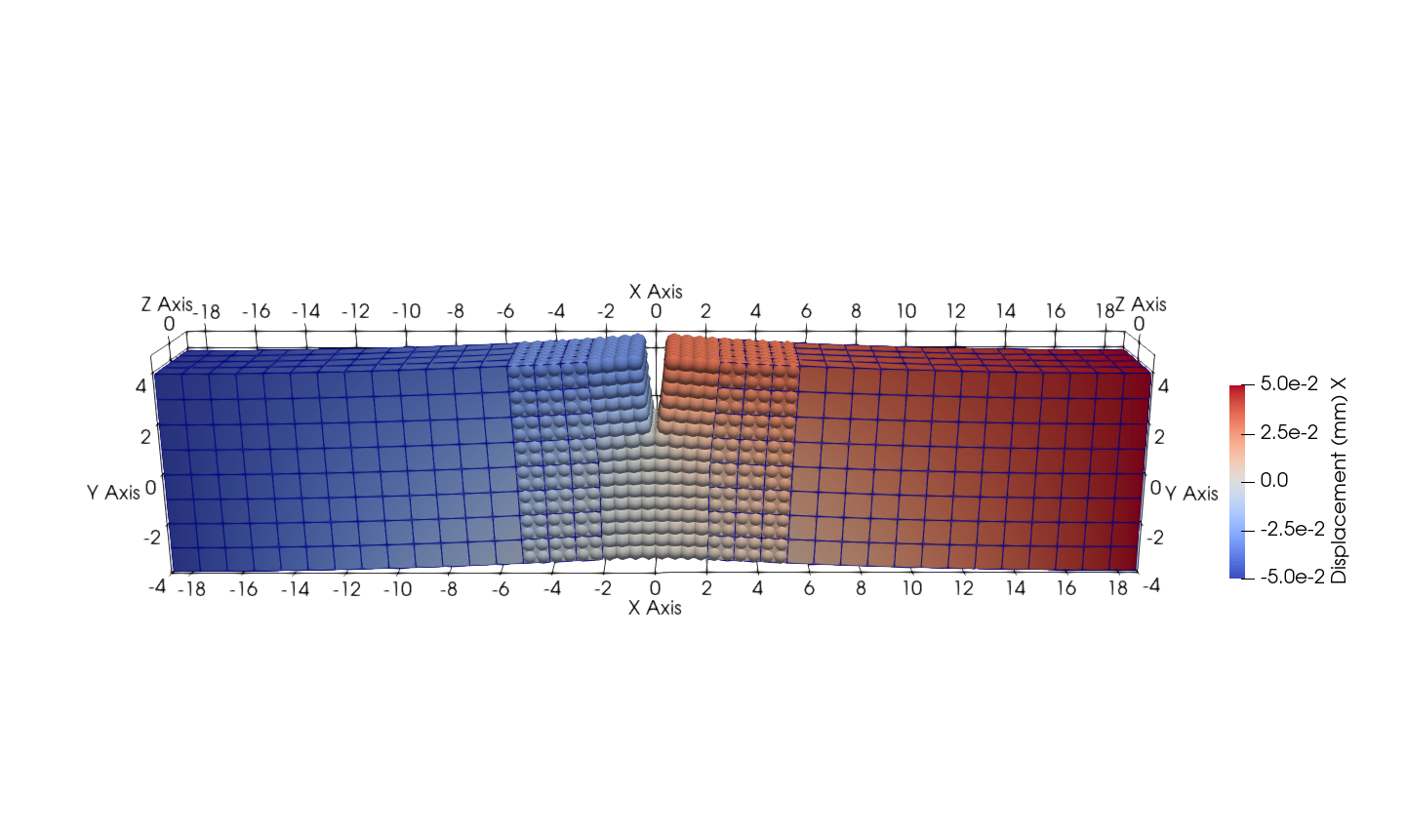}
     \subcaption{$X$ component displacements.}\label{fig::DirichletDirichletBarMagnitude}  
   }   
   \end{subfigure}   

  \begin{subfigure}{\textwidth}
    {   \centering
     \includegraphics[scale=0.35,trim=3.1cm 8.1cm 2.1cm 10.3cm,clip]{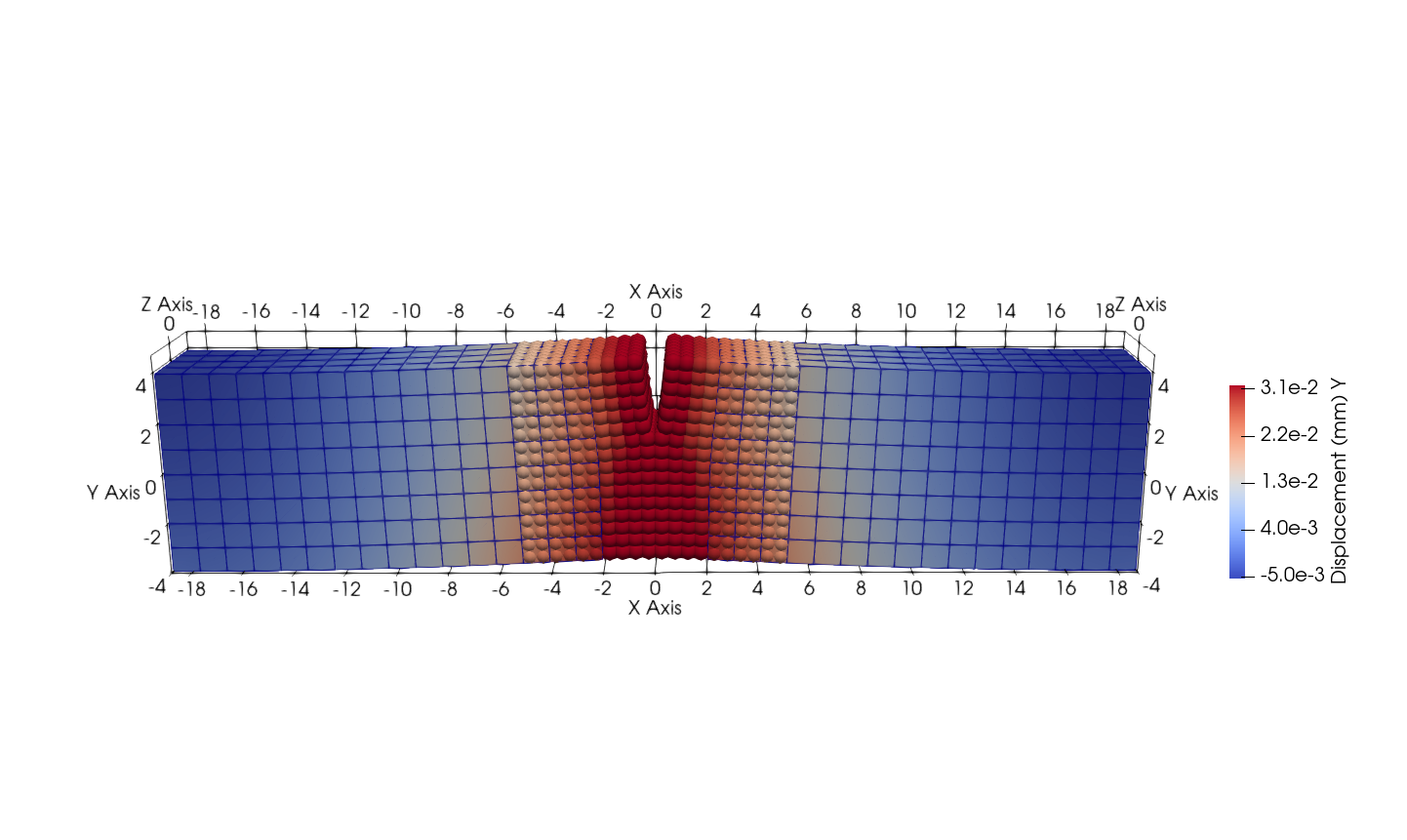}
     \subcaption{$Y$ component displacements}\label{fig::DirichletDirichletBarX}
     }
     \end{subfigure}

  \begin{subfigure}{\textwidth}
    {   \centering
     \includegraphics[scale=0.35,trim=3.1cm 8.1cm 2.1cm 10.3cm,clip]{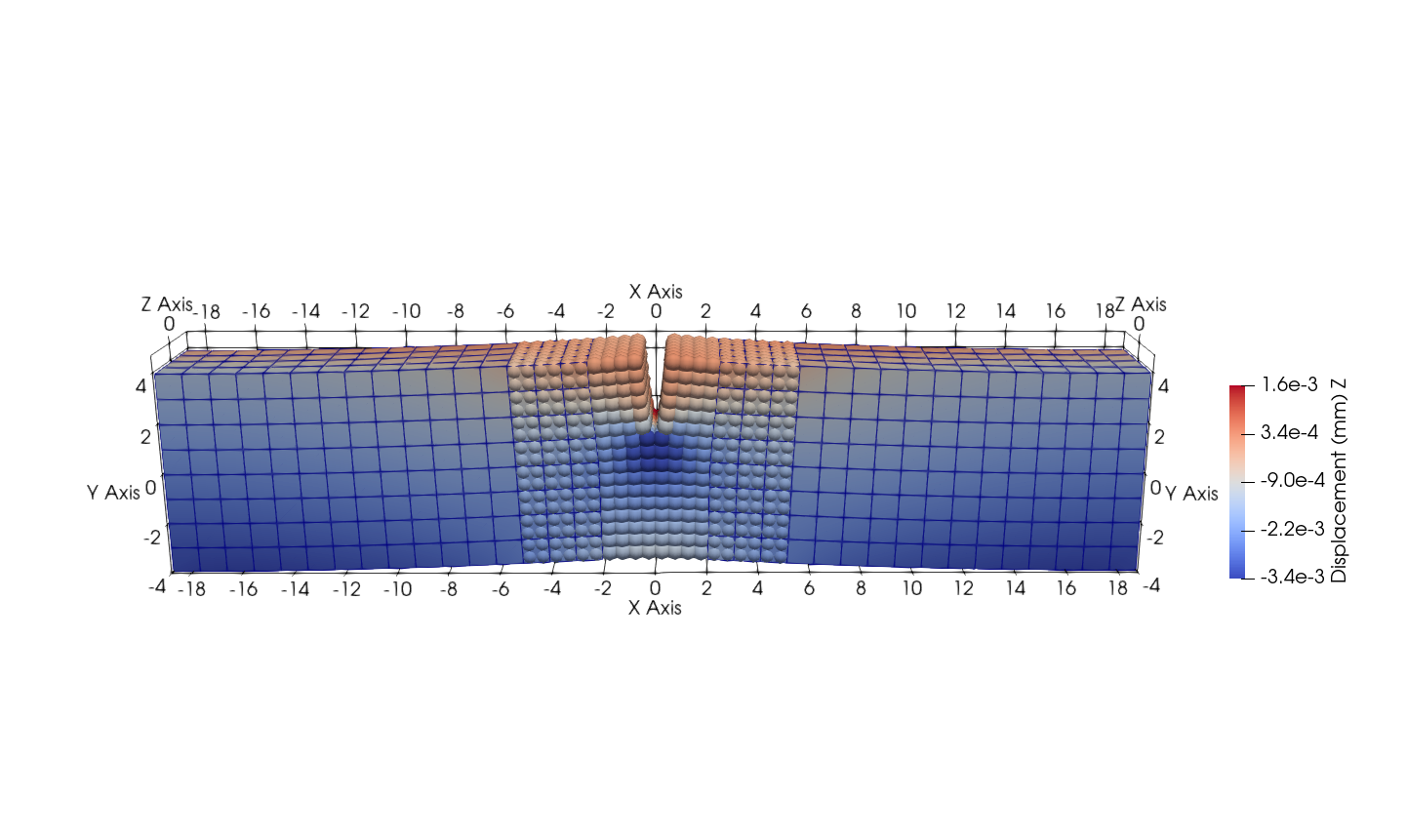}
     \subcaption{$Z$ component displacements}\label{fig::DirichletDirichletBarY}
     }
     \end{subfigure}

    \end{center}
    \caption{Results of the bar experiment with prescribed displacements of \SI{-0.05}{\milli\meter} at $x=\SI{-18.0}{\milli\meter}$ and \SI{0.05}{\milli\meter} at $x=\SI{18.0}{\milli\meter}$ (see Figure \ref{fig:DirichletBarIllustration}). Displacements are magnified by a factor of $15$ for the purpose of illustration.}\label{fig::dirichletdirichletexperiment}
\end{figure}

This experiment considers a prenotched bar described by the geometry $\Omega := [-18,18] \times [-4,4] \times [-2,2]$, where the unit of length is \SI{}{\milli\meter}. See Figure \ref{fig:DirichletBarIllustration} for an illustration of a cross-section along the $XY$ plane. The local solution domain, boundary regions, and control regions are given by
\begin{equation}
\begin{split}
   \Omega_l :={}& \left((-18,-2) \times [-4,4] \times [-2,2] \right) \cup \left( (2,18) \times [-4,4] \times [-2,2] \right) \\
   \Gamma_D :={}& \left(\left\{-18\right\} \times [-4,4] \times [-2,2]\right) \cup \left(\left\{18\right\} \times [-4,4] \times [-2,2] \right) \\
   \Gamma_c :={}& \left(\left\{-2\right\} \times [-4,4] \times [-2,2]\right) \cup \left(\left\{2\right\} \times [-4,4] \times [-2,2] \right)
\end{split}
\end{equation}
The nonlocal solution domain and control regions are given by
\begin{equation}
    \begin{split}
        \omega_n :={}& (-3.5,3.5) \times [-4,4] \times [-2,2] \\
        \eta_c :={}& \left([-5,-3.5] \times [-4,4] \times [-2,2] \right) \cup \left([3.5,5] \times [-4,4] \times [-2,2] \right)
    \end{split}
\end{equation}
In addition, a prenotch region $P$ is described by
\begin{equation}
    P:= \left\{0 \right\} \times [1,4] \times [-2,2].
\end{equation}
The prenotch is modeled in the numerical experiment by omitting all bonds crossing $P$. Our choice of domain decomposition is dictated by the fact that the nonlocal model is advantageous in the vicinity of the prenotch, while the local model is sufficient far from it.  The nonlocal and local computational domains are further illustrated in Figure \ref{fig::dirichletdirichletexperiment}, where material points in the nonlocal model are rendered as spheres and the hexahedral FE mesh corresponds to the local model.

A prescribed displacement boundary condition is applied in the $x$ direction on the local boundary region $\Gamma_D$, specifically \SI{-0.05}{\milli\meter} at $x = \SI{-18.0}{\milli\meter}$ and \SI{0.05}{\milli\meter} at $x = \SI{18.0}{\milli\meter}$, resulting in tensile loading. Stress-free conditions are employed on the remaining surfaces, i.e., a homogeneous Neumann volume constraint and boundary condition for the nonlocal and local problems, respectively. To eliminate rigid body modes, additional zero displacement boundary conditions are applied in the $y$ direction along the edges defined by $x = \SI{-18.0}{\milli\meter}$ , $y = \SI{-4.0}{\milli\meter}$ and $x = \SI{18.0}{\milli\meter}$, $y = \SI{-4.0}{\milli\meter}$, and in the $z$ direction along the edges defined by $x = \SI{-18.0}{\milli\meter}$, $z = \SI{-2.0}{\milli\meter}$ and $x = \SI{18.0}{\milli\meter}$, $z = \SI{-2.0}{\milli\meter}$. The results of the experiment are shown in Figure \ref{fig::dirichletdirichletexperiment}. As expected, the influence of the crack on the displacement is restricted predominantly to the nonlocal region, and the OB coupling provides a smooth transition between the local and nonlocal models.

\subsubsection{Compact tension experiment}
To investigate the viability of the OB coupling method on a more realistic engineering geometry, we consider the compact tension test specimen illustrated in Figure~\ref{fig:CompactTensionIllustration}.  The compact tension test is a common laboratory experiment for the evaluation of material properties, such as fracture toughness, that are related to material failure.  We employ a compact tension test geometry based on the ASTM E399 - 20 standard, although we reduce the thickness to decrease computational expense. As in the case of the prenotched bar in Section~\ref{sec:prenotched_bar}, the compact tension test simulation demonstrates the ability of the coupling strategy to permit models in which peridynamics is applied to regions where cracks are present, and a less computationally expensive classical FE model is applied elsewhere.  The local-nonlocal interface is placed sufficiently far from the crack tip such that the material response in the overlap region is smooth and well behaved.  Prescribed displacement boundary conditions are applied at the holes, which is representative of the physical laboratory experiment in which pins connected to a material testing machine are used to load the specimen.  It is noted that in this configuration, loading is transmitted to the peridynamic region only though the coupling interface; nonzero user-prescribed boundary conditions are confined to the local model.  The results of the simulation are presented in Figure \ref{Fig::CompactTensionExperiment}.  The deformations are consistent with the applied loading, and no unphysical artifacts are apparent in the overlap region.

The simulations carried out in the present study are static, i.e., the coupled \emph{Peridigm}-\emph{Albany} code solves for the equilibrium configuration corresponding to a single set of prescribed boundary conditions.  Extending the code capabilities to enable quasi-static simulations that include stable crack propagation emanating from the prenotch is a subject of future work.

\begin{figure}
\begin{center}
\begin{tikzpicture}[scale=0.7]

\def\OmegaLColor{rgb,255:red,60; green,108; blue,199};
\def\omeganColor{rgb,255:red,162; green,214; blue,235};
\def\etacColor{yellow};
\def\etaDColor{pink};
\def\GammacColor{black!10!red};
\def\GammaDColor{brown};
\def\prenotchColor{white};

\def\omegalregion{(-6.25cm,-6.0cm) rectangle (6.25cm,6.0cm)}; 
\def\omegaletacregion{(-2.7cm,-1.8cm) rectangle (0.9cm,1.8cm)};
\def\omegalomeganregion{(-2.4cm,-1.5cm) rectangle (0.6cm,1.5cm)};
\def\omeganregion{(-2.1cm,-1.2cm) rectangle (0.3cm,1.2cm)};
\def\prenotchregion{(6.35,0.25) -- (1.25,0.25) -- (0.05,0.0) -- (1.25,-0.25) -- (6.35,-0.25) -- cycle};
\def\topcircleregion{(3.75,2.75) circle (1.25cm);};
\def\botcircleregion{(3.75,-2.75) circle (1.25cm);};

\fill [fill=\OmegaLColor] \omegalregion;

\begin{scope} 
  \clip \omegaletacregion;

  \foreach \x in {-100,...,20}
    {\pgfmathsetmacro{\a}{-5+2*\x/4}
     \pgfmathsetmacro{\b}{-5/4-5/4*\a}
     \pgfmathsetmacro{\c}{-5+2*\x/4+1/4}
     \pgfmathsetmacro{\d}{-5/4-5/4*\c}
     \fill[fill={\OmegaLColor}] (\a,-30) -- (30,\b) -- (30,-30) -- cycle;
     \fill[fill={\etacColor}] (\c,-30) -- (30,\d) -- (30,-30) -- cycle;}
\end{scope}

\begin{scope} 
  \clip \omegalomeganregion;

  \foreach \x in {-100,...,30}
    {\pgfmathsetmacro{\a}{-5+2*\x/4}
     \pgfmathsetmacro{\b}{-5/4-5/4*\a}
     \pgfmathsetmacro{\c}{-5+2*\x/4+1/4}
     \pgfmathsetmacro{\d}{-5/4-5/4*\c}
     \fill[fill={\OmegaLColor}] (\a,-30) -- (30,\b) -- (30,-30) -- cycle;
     \fill[fill={\omeganColor}] (\c,-30) -- (30,\d) -- (30,-30) -- cycle;}
\end{scope}

\fill [fill=\omeganColor] \omeganregion;


\draw[color={black}, line width=0.1mm] (-2.7,-1.8cm) -- (-2.7,1.8cm);
\draw[color={black}, line width=0.1mm] (0.9,-1.8cm) -- (0.9,1.8cm);
\draw[color={black}, line width=0.1mm] (-2.4,-1.5cm) -- (-2.4,1.5cm);
\draw[color={black}, line width=0.1mm] (0.6,-1.5cm) -- (0.6,1.5cm);
\draw[color={black}, line width=0.1mm] (-2.7,1.8cm) -- (0.9,1.8cm);
\draw[color={black}, line width=0.1mm] (-2.7,-1.8cm) -- (0.9,-1.8cm);
\draw[color={black}, line width=0.1mm] (-2.4,1.5cm) -- (0.6,1.5cm);
\draw[color={black}, line width=0.1mm] (-2.4,-1.5cm) -- (0.6,-1.5cm);

\draw[color={\GammacColor}, line width=0.25mm] (-2.1,-1.2cm) -- (-2.1,1.2cm);
\draw[color={\GammacColor}, line width=0.25mm] (0.3,-1.2cm) -- (0.3,1.2cm);
\draw[color={\GammacColor}, line width=0.25mm] (-2.1,-1.2cm) -- (0.3,-1.2cm);
\draw[color={\GammacColor}, line width=0.25mm] (-2.1,1.2cm) -- (0.3,1.2cm);

\draw[color={black}, line width=0.25mm] (-6.25,-6.0cm) -- (-6.25,6.0cm);
\draw[color={black}, line width=0.25mm] (6.25,-6.0cm) -- (6.25,6.0cm);
\draw[color={black}, line width=0.25mm] (-6.25,-6.0cm) -- (6.25,-6.0cm);
\draw[color={black}, line width=0.25mm] (-6.25,6.0cm) -- (6.25,6.0cm);

\fill [fill=\prenotchColor] \prenotchregion;
\draw[color={black}, line width=0.1mm] (6.25,0.25) -- (1.25,0.25) -- (0.05,0.0) -- (1.25,-0.25) -- (6.25,-0.25);

\fill [fill=\prenotchColor] \topcircleregion;
\fill [fill=\prenotchColor] \botcircleregion;
\draw[color={black}, line width=0.25mm] \topcircleregion;
\draw[color={black}, line width=0.25mm] \botcircleregion;

\begin{scope} 
  \clip (2.75,2.75) rectangle (4.75,8.0);
\draw[color={\GammaDColor}, line width=0.25mm] \topcircleregion;
\end{scope}
\begin{scope} 
  \clip (2.75,-2.75) rectangle (4.75,-8.0);
\draw[color={\GammaDColor}, line width=0.25mm] \botcircleregion;
\end{scope}

\draw[<->] (-8.5,-6) -- (-8.5,6) node[midway, left, rotate=90,yshift=1.2ex,xshift=4.0ex]{\footnotesize \SI{12.0}{\milli\meter}};
\draw[dotted] (-8.5,-6) -- (-6.25,-6);
\draw[dotted] (-8.5,6) -- (-6.25,6);

\draw[<->] (-6.55,-1.2) -- (-6.55,1.2) node[midway, left, rotate=90,yshift=1.2ex,xshift=4.0ex]{\footnotesize \SI{2.4}{\milli\meter}};
\draw[dotted] (-6.55,-1.2) -- (-2.1,-1.2);
\draw[dotted] (-6.55,1.2) -- (-2.1,1.2);

\draw[<->] (-7.2,-1.5) -- (-7.2,1.5) node[midway, left, rotate=90,yshift=1.2ex,xshift=4.0ex]{\footnotesize \SI{3.0}{\milli\meter}};
\draw[dotted] (-7.2,-1.5) -- (-2.4,-1.5);
\draw[dotted] (-7.2,1.5) -- (-2.4,1.5);

\draw[<->] (-7.85,-1.8) -- (-7.85,1.8) node[midway, left, rotate=90,yshift=1.2ex,xshift=4.0ex]{\footnotesize \SI{3.6}{\milli\meter}};
\draw[dotted] (-7.85,-1.8) -- (-2.7,-1.8);
\draw[dotted] (-7.85,1.8) -- (-2.7,1.8);

\draw[<->] (-6.25,7.7) -- (6.25,7.7) node[midway, above,yshift=-0.3ex]{\footnotesize \SI{12.5}{\milli\meter}};
\draw[dotted] (-6.25,7.7) -- (-6.25,6.0);
\draw[dotted] (6.25,7.7) -- (6.25,6.0);

\draw[<->] (2.5,6.4) -- (6.25,6.4) node[midway, above,yshift=-0.3ex]{\footnotesize \SI{3.75}{\milli\meter}};
\draw[dotted] (2.5,6.4) -- (2.5,2.75);

\draw[<->] (0.05,7.05) -- (6.25,7.05) node[midway, above,yshift=-0.3ex]{\footnotesize \SI{6.2}{\milli\meter}};
\draw[dotted] (0.05,7.05) -- (0.05,0.0);

\draw[<->] (6.75,-0.25) -- (6.75,0.25) node[midway, left, rotate=-90,yshift=1.5ex,xshift=4.0ex]{\footnotesize \SI{0.5}{\milli\meter}};
\draw[dotted] (6.25,-0.25) -- (6.75,-0.25);
\draw[dotted] (6.25,0.25) -- (6.75,0.25);

\draw[<->] (6.75,1.5) -- (6.75,4.0) node[midway, left, rotate=-90,yshift=1.2ex,xshift=4.0ex]{\footnotesize \SI{2.5}{\milli\meter}};
\draw[dotted] (3.75,1.5) -- (6.75,1.5);
\draw[dotted] (3.75,4.0) -- (6.75,4.0);

\draw[<->] (6.75,-6.0) -- (6.75,-2.75) node[midway, left, rotate=-90,yshift=1.2ex,xshift=4.0ex]{\footnotesize \SI{3.25}{\milli\meter}};
\draw[dotted] (6.0,-6.0) -- (6.75,-6.0);
\draw[dotted] (5.0,-2.75) -- (6.75,-2.75);

  \draw[color={\GammaDColor}, line width=0.5mm] (9.0,4.25) -- (10.0,4.25) node[right,text=black] {$\Gamma_D$};
  \draw[color={\GammacColor}, line width=0.5mm] (9.0,3.25) -- (10.0,3.25) node[right,text=black] {$\Gamma_c$};
  \draw[color={\OmegaLColor},line width=2.5mm] (9.0,2.25) -- (10.0,2.25) node[right,text=black] {$\Omega_l$};

  \draw[color={\omeganColor},line width=2.5mm] (9.0,-2.25) -- (10.0,-2.25) node[right,text=black] {$\omega_n$};
  \draw[color={\etacColor},line width=2.5mm] (9.0,-3.25) -- (10.0,-3.25) node[right,text=black] {$\eta_c$};

\end{tikzpicture}
\caption{$XY$ cross section of the geometry for the compact tension test. Dirichlet boundary conditions are imposed on $\Gamma_D$, where we impose a displacement of $\langle 0.0,0.03,0.0 \rangle\,\SI{}{\milli\meter}$ in the top hole and a displacement of $\langle 0.0,-0.03,0.0 \rangle\,\SI{}{\milli\meter}$ in the bottom hole. }\label{fig:CompactTensionIllustration}
\end{center}
\end{figure}
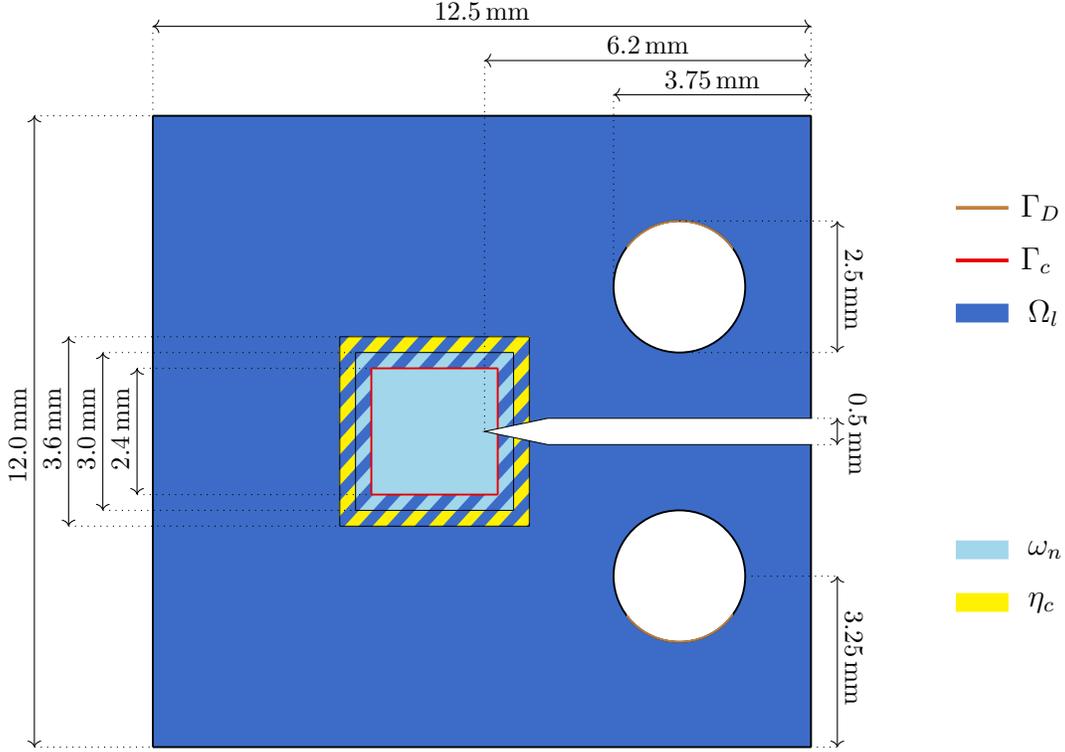
\begin{figure}
    \begin{center}
  \begin{subfigure}{0.48\textwidth}
    {   \centering
     \includegraphics[trim=0.7cm 0.7cm 4.8cm 2.1cm,clip,width=1.0\textwidth]{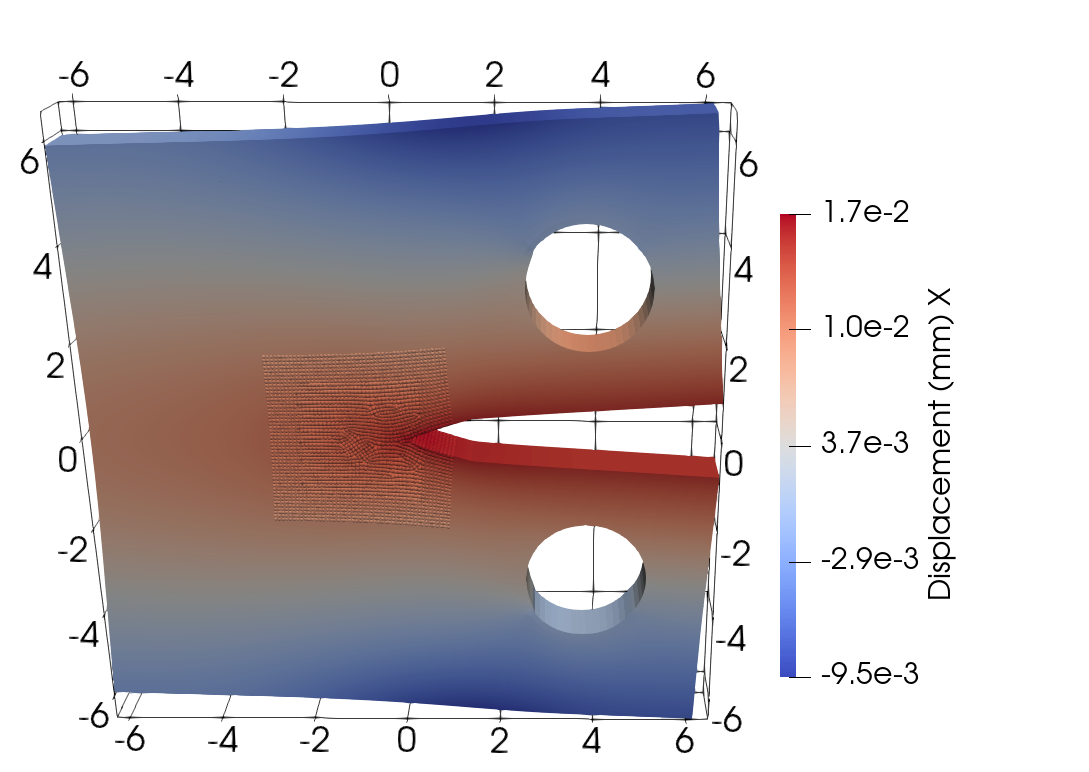}
     \subcaption{Displacement in $x$.}\label{fig::CompactTensionX}
     }
     \end{subfigure}
\hspace*{5pt}  \begin{subfigure}{0.48\textwidth}
    {   \centering
     \includegraphics[trim=0.7cm 0.7cm 4.8cm 2.1cm,clip,width=1.0\textwidth]{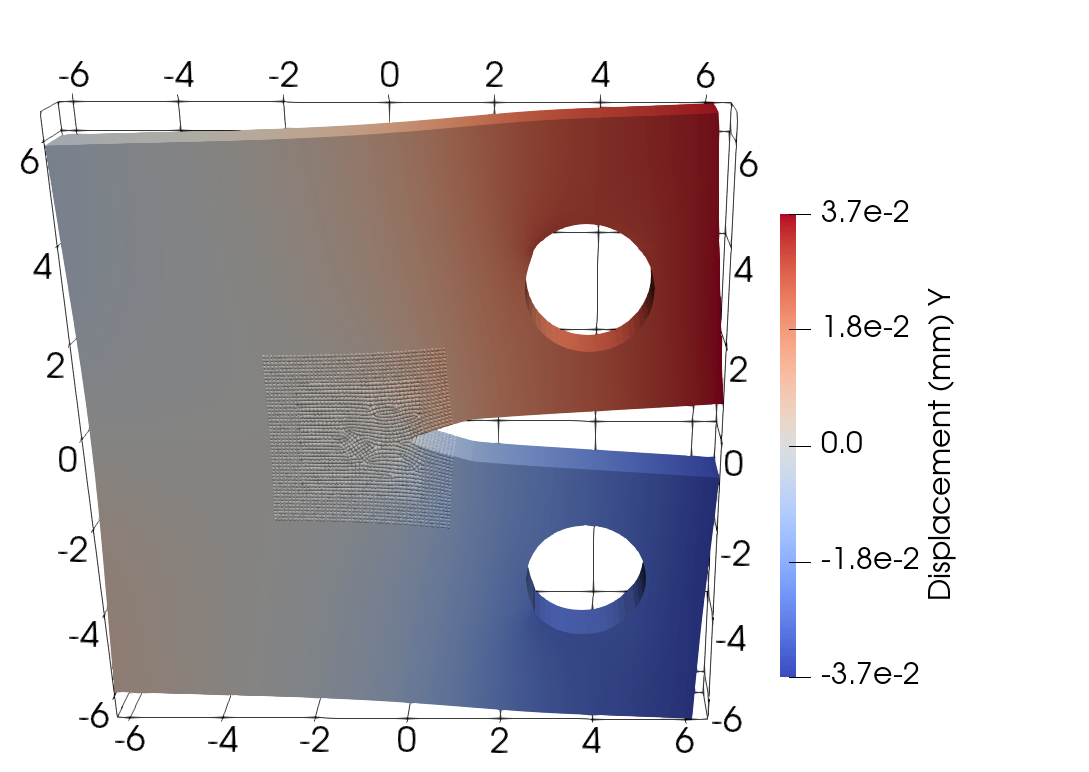}
     \subcaption{Displacement in $y$.}\label{fig::CompactTensionY}
     }
     \end{subfigure}

  \begin{subfigure}{0.48\textwidth}
    {   \centering
     \includegraphics[trim=0.7cm 0.7cm 4.8cm 2.1cm,clip,width=1.0\textwidth]{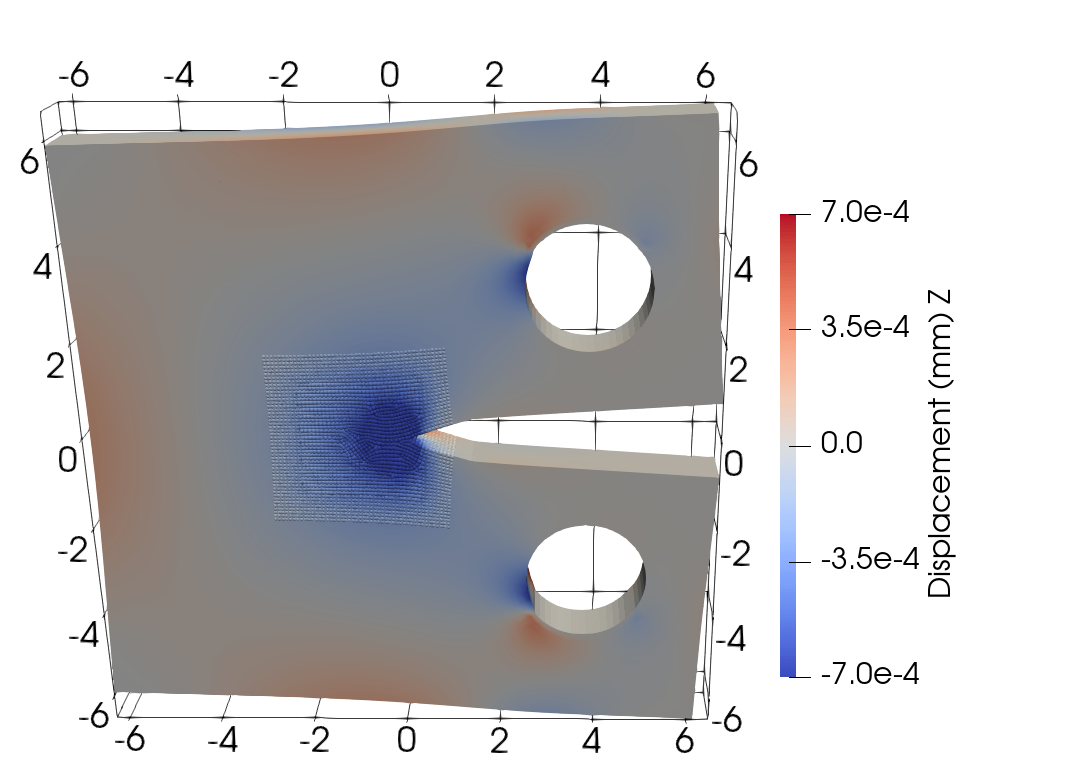}
     \subcaption{Displacement in $z$.}\label{fig::CompactTensionZ}
     }
     \end{subfigure}
\hspace*{5pt}  \begin{subfigure}{0.48\textwidth}
    {   \centering
     \includegraphics[trim=0.7cm 0.7cm 4.8cm 2.1cm,clip,width=1.0\textwidth]{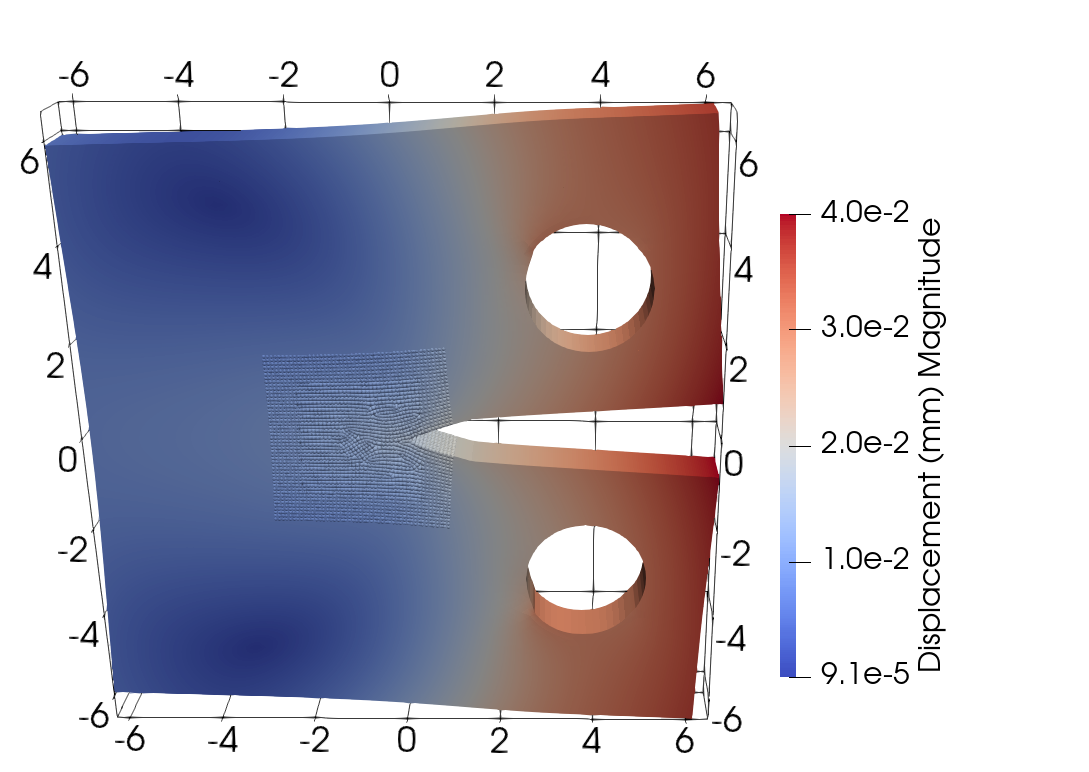}
     \subcaption{Magnitude of displacement.}\label{fig::CompactTensionMag}
     }
     \end{subfigure}

    \end{center}
    \caption{Compact Tension Experiment with Dirichlet Boundary Conditions}\label{Fig::CompactTensionExperiment}
\end{figure}

\subsubsection{Conversion of local to nonlocal boundary conditions in a prenotched bar}
A notable use case for LtN coupling is the conversion of local boundary conditions into nonlocal volume constraints.  This is important because surface traction boundary conditions, commonly applied in traditional FE models to represent pressure loading, or similar, cannot be applied directly to a nonlocal model.  Instead, to ensure uniqueness of solution for the nonlocal problem, tractions must be converted from two-dimensional surface loads for use in a local model to three-dimensional volume constraints for use in a nonlocal model.

\begin{figure}
\begin{center}
\begin{tikzpicture}[scale=0.25]

\def\OmegaLColor{rgb,255:red,60; green,108; blue,199};
\def\omeganColor{rgb,255:red,162; green,214; blue,235};
\def\etacColor{yellow};
\def\etaDColor{pink};
\def\GammacColor{black!10!red};
\def\GammaDColor{brown};

\def\omeganregion{(-8.0cm,-4.0cm) rectangle (8.0cm,4.0cm)};
\def\omegalleftregion{(-16cm,-4.0cm) rectangle (-8.0cm,4.0cm)}; 
\def\omegalrightregion{(8.0cm,-4.0cm) rectangle (16.0cm,4.0cm)}; 
\def\omegalomeganregionl{(-14.5cm,-4.0cm) rectangle (-8.0cm,4.0cm)};
\def\omegalomeganregionr{(8.0cm,-4.0cm) rectangle (14.5cm,4.0cm)};
\def\omegaletacregionl{(-16.0cm,-4.0cm) rectangle (-14.5cm,4.0cm)};
\def\omegaletacregionr{(14.5cm,-4.0cm) rectangle (16.0cm,4.0cm)};

\fill [fill=\omeganColor] \omeganregion;

\fill [fill=\OmegaLColor] \omegalleftregion;
\fill [fill=\OmegaLColor] \omegalrightregion;


\begin{scope} 
  \clip \omegalomeganregionl;

  \foreach \x in {-100,...,30}
    {\pgfmathsetmacro{\a}{-5+2*\x/4}
     \pgfmathsetmacro{\b}{-5/4-5/4*\a}
     \pgfmathsetmacro{\c}{-5+2*\x/4+1/4}
     \pgfmathsetmacro{\d}{-5/4-5/4*\c}
     \fill[fill={\OmegaLColor}] (\a,-30) -- (30,\b) -- (30,-30) -- cycle;
     \fill[fill={\omeganColor}] (\c,-30) -- (30,\d) -- (30,-30) -- cycle;}
\end{scope}
\begin{scope} 
  \clip \omegalomeganregionr;

  \foreach \x in {-60,...,60}
    {\pgfmathsetmacro{\a}{-5+2*\x/4}
     \pgfmathsetmacro{\b}{-5/4-5/4*\a}
     \pgfmathsetmacro{\c}{-5+2*\x/4+1/4}
     \pgfmathsetmacro{\d}{-5/4-5/4*\c}
     \fill[fill={\OmegaLColor}] (\a,-30) -- (30,\b) -- (30,-30) -- cycle;
     \fill[fill={\omeganColor}] (\c,-30) -- (30,\d) -- (30,-30) -- cycle;}
\end{scope}

\begin{scope} 
  \clip \omegaletacregionl;

  \foreach \x in {-100,...,20}
    {\pgfmathsetmacro{\a}{-5+2*\x/4}
     \pgfmathsetmacro{\b}{-5/4-5/4*\a}
     \pgfmathsetmacro{\c}{-5+2*\x/4+1/4}
     \pgfmathsetmacro{\d}{-5/4-5/4*\c}
     \fill[fill={\OmegaLColor}] (\a,-30) -- (30,\b) -- (30,-30) -- cycle;
     \fill[fill={\etacColor}] (\c,-30) -- (30,\d) -- (30,-30) -- cycle;}
\end{scope}
\begin{scope} 
  \clip \omegaletacregionr;

  \foreach \x in {-40,...,30}
    {\pgfmathsetmacro{\a}{-5+2*\x/4}
     \pgfmathsetmacro{\b}{-5/4-5/4*\a}
     \pgfmathsetmacro{\c}{-5+2*\x/4+1/4}
     \pgfmathsetmacro{\d}{-5/4-5/4*\c}
     \fill[fill={\OmegaLColor}] (\a,-30) -- (30,\b) -- (30,-30) -- cycle;
     \fill[fill={\etacColor}] (\c,-30) -- (30,\d) -- (30,-30) -- cycle;}
\end{scope}

\draw[color={black}, line width=0.1mm] (-14.5,-4.0cm) -- (-14.5,4.0cm);
\draw[color={black}, line width=0.1mm] (14.5,-4.0cm) -- (14.5,4.0cm);
\draw[color={black}, line width=0.1mm] (-8.0,-4.0cm) -- (-8.0,4.0cm);
\draw[color={black}, line width=0.1mm] (8.0,-4.0cm) -- (8.0,4.0cm);
\draw[color={black}, line width=0.1mm] (-16.0,-4.0cm) -- (16.0,-4.0cm);
\draw[color={black}, line width=0.1mm] (-16.0,4.0cm) -- (16.0,4.0cm);

\draw[color={white}, line width=0.25mm] (0.0,1.0cm) -- (0.0,4.02cm);

\draw[color={\GammaDColor}, line width=0.5mm] (-16.0,-4.02cm) -- (-16.0,4.02cm);
\draw[color={\GammaDColor}, line width=0.5mm] (16.0,-4.02cm) -- (16.0,4.02cm);

\draw[color={\GammacColor}, line width=0.5mm] (-8.0,-4.02cm) -- (-8.0,4.02cm);
\draw[color={\GammacColor}, line width=0.5mm] (8.0,-4.02cm) -- (8.0,4.02cm);

 \draw[color={\GammaDColor}, line width=0.5mm] (19.0,2.25) -- (21.0,2.25) node[right,text=black] {$\Gamma_D$};
 \draw[color={\GammacColor}, line width=0.5mm] (19.0,0.0) -- (21.0,0.0) node[right,text=black] {$\Gamma_c$};
 \draw[color={\OmegaLColor},line width=2.5mm] (19.0,-2.25) -- (21.0,-2.25) node[right,text=black] {$\Omega_l$};

 \draw[color={\omeganColor},line width=2.5mm] (-21.0,2.25) -- (-19.0,2.25) node[left,text=black,xshift=-2ex] {$\omega_n$};
 \draw[color={\etacColor},line width=2.5mm] (-21.0,-2.25) -- (-19.0,-2.25) node[left,text=black,xshift=-2ex] {$\eta_c$};

\draw[<->] (14.5,5.0) -- (16.0,5.0) node[midway, above, rotate=0,yshift=-0.3ex,xshift=-2.4ex]{\footnotesize $1.5$mm};
\draw[dotted] (14.5,4.0) -- (14.5,5.0);

\draw[<->] (8.0,6.5) -- (16.0,6.5) node[midway, above, rotate=0,yshift=-0.3ex,xshift=0.0ex]{\footnotesize $8.0$mm};
\draw[dotted] (8.0,4.0) -- (8.0,6.5);

\draw[<->] (0,8.0) -- (16.0,8.0) node[midway, above, rotate=0,yshift=-0.3ex,xshift=0.0ex]{\footnotesize $16.0$mm};
\draw[dotted] (0.0,4.0) -- (0.0,8.0);

\draw[<->] (-16.0,9.5) -- (16.0,9.5) node[midway, above, rotate=0,yshift=-0.3ex,xshift=0.0ex]{\footnotesize $32.0$mm};
\draw[dotted] (-16.0,4.0) -- (-16.0,9.5);
\draw[dotted] (16.0,4.0) -- (16.0,9.5);

\end{tikzpicture}
\caption{$XY$ cross section geometry for the prenotched bar experiment in which a Dirichlet boundary condition is imposed on the right end of the bar and a Neumann boundary condition is imposed on the left side of the bar.}\label{fig:NeumannBarIllustration}
\end{center}
\end{figure}
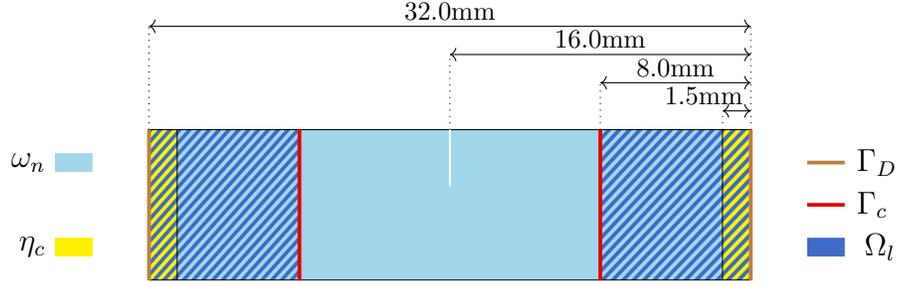

We demonstrate the OB coupling method for conversion of traction boundary conditions to nonlocal volume constraints using the model of a bar illustrated in Figure~\ref{fig:NeumannBarIllustration}.  The computational domain is defined as  $\Omega := [-16,16] \times [-4,4] \times [-2,2]$, where the unit of length is millimeters. The local solution domain, boundary regions, and control regions are, respectively, given by
\begin{equation}
\begin{split}
   \Omega_l :={}& \left((-16,-8) \times [-4,4] \times [-2,2] \right) \cup \left( (8,16) \times [-4,4] \times [-2,2] \right) \\
   \Gamma_D :={}& \left(\left\{-16\right\} \times [-4,4] \times [-2,2]\right) \cup \left(\left\{16\right\} \times [-4,4] \times [-2,2] \right) \\
   \Gamma_c :={}& \left(\left\{-8\right\} \times [-4,4] \times [-2,2]\right) \cup \left(\left\{8\right\} \times [-4,4] \times [-2,2] \right)
\end{split}
\end{equation}
The nonlocal solution domain and control regions are, respectively, given by
\begin{equation}
    \begin{split}
        \omega_n :={}& (-14.5,14.5) \times [-4,4] \times [-2,2] \\
        \eta_c :={}& \left([-16,-14.5] \times [-4,4] \times [-2,2] \right) \cup \left([14.5,16] \times [-4,4] \times [-2,2] \right)
    \end{split}
\end{equation}
In addition, a prenotch region $P$ is described by
\begin{equation}
    P:= \left\{0 \right\} \times [1,4] \times [-2,2].
\end{equation}
In the numerical experiment, all bonds crossing $P$ are omitted from the simulation. On the local boundary region $\Gamma_D$, a prescribed traction of $\tau =$ \SI{-1700.0}{\mega\pascal} in the $x$ direction at $x=$ \SI{-16.0}{\milli\meter} and a prescribed displacement of \SI{0.0}{\milli\meter} in the $x$ direction at $x=$ \SI{16.0}{\milli\meter} are imposed. To eliminate rigid body modes, additional zero displacement boundary conditions are applied in the $y$ direction along the edge defined by $x=$ \SI{16.0}{\milli\meter}, $y=$ \SI{-4.0}{\milli\meter}, and in the $z$ direction along the edge defined by $x=$ \SI{16.0}{\milli\meter}, $z=$ \SI{-2.0}{\milli\meter}.

In contrast to the prenotched bar experiment in Section~\ref{sec:prenotched_bar}, the nonlocal region is not restricted to the vicinity of the prenotch and instead encompasses the entire bar. This approach allows us to observe the specific nonlocal volume constrains that correspond to the local boundary conditions at the ends of the bar.  Of particular interest are the nonlocal volume constraints that correspond to the traction boundary condition applied to the face at $x=$ \SI{-16.0}{\milli\meter}.  The nonlocal volume constraints, which are virtual constraints determined by solution of the OB coupling problem, give insight into how nonlocal volume constrains could be applied in a purely nonlocal simulation to reproduce the effect of traction loading in a comparable local simulation. This conversion technique is similar to the approaches introduced in \cite{DEliaNeumann2019,DElia2021prescription}.

\begin{figure}
    \begin{center}
  \begin{subfigure}{\textwidth}
    {   \centering
     \includegraphics[scale=0.35,trim=0.7cm 6.5cm 5.3cm 9.5cm,clip]{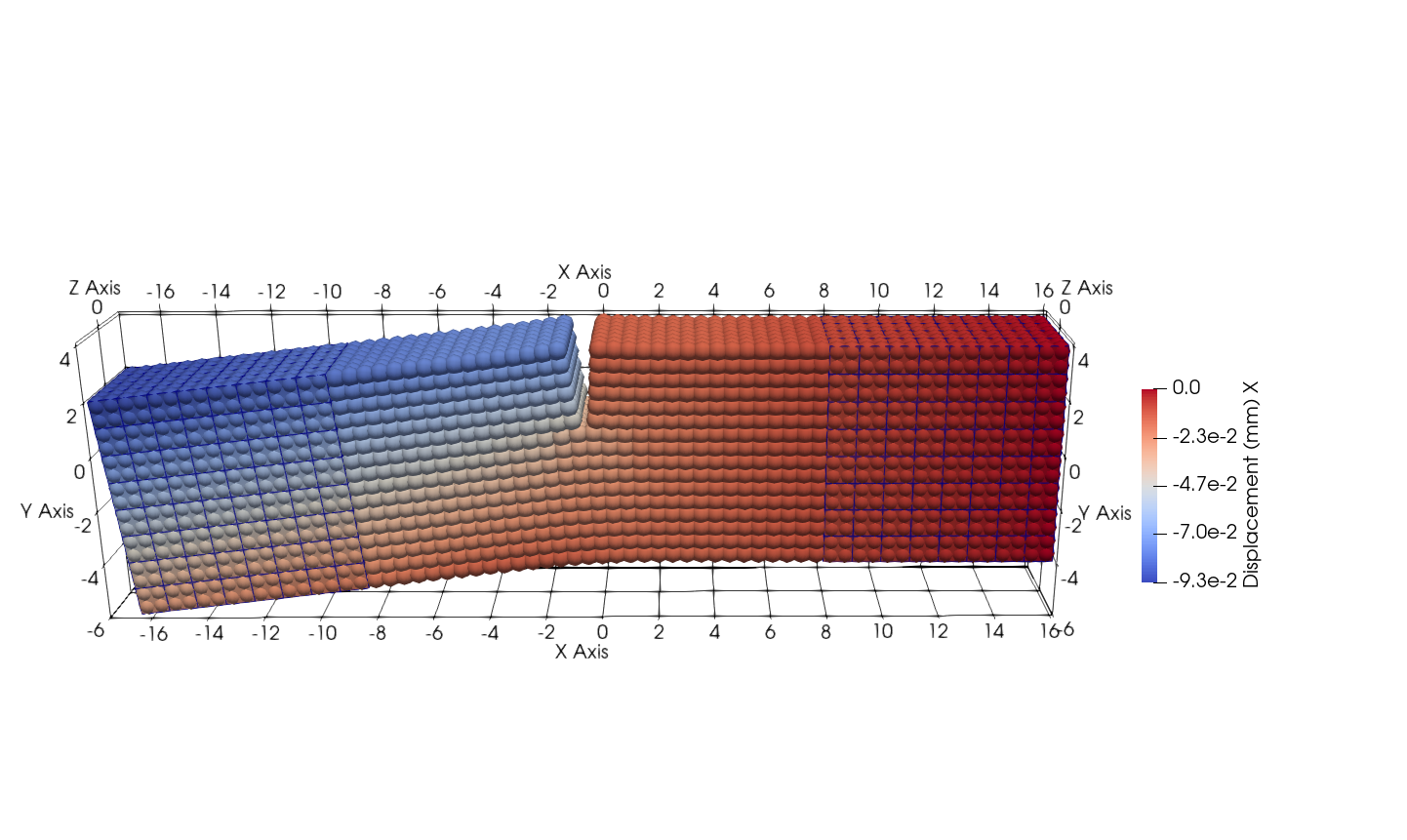}
     \subcaption{Displacement in $x$}\label{fig::DirichletNeumannBarX}
     }
     \end{subfigure}

  \begin{subfigure}{\textwidth}
    {   \centering
     \includegraphics[scale=0.35,trim=0.7cm 6.5cm 5.3cm 9.5cm,clip]{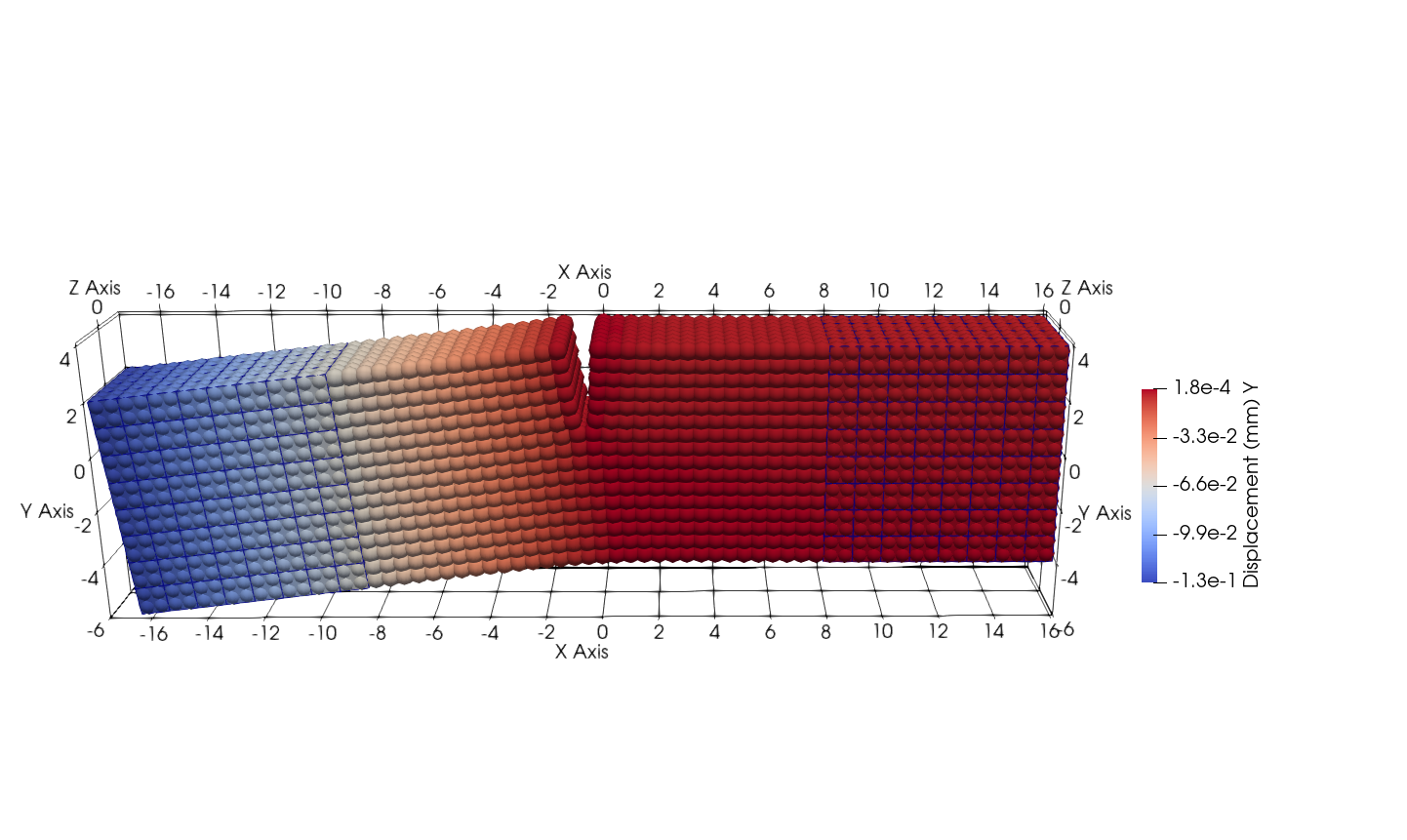}
     \subcaption{Displacement in $y$}\label{fig::DirichletNeumannBarY}
     }
     \end{subfigure}
     
  \begin{subfigure}{\textwidth}
    {   \centering
     \includegraphics[scale=0.35,trim=0.7cm 6.5cm 5.3cm 9.5cm,clip]{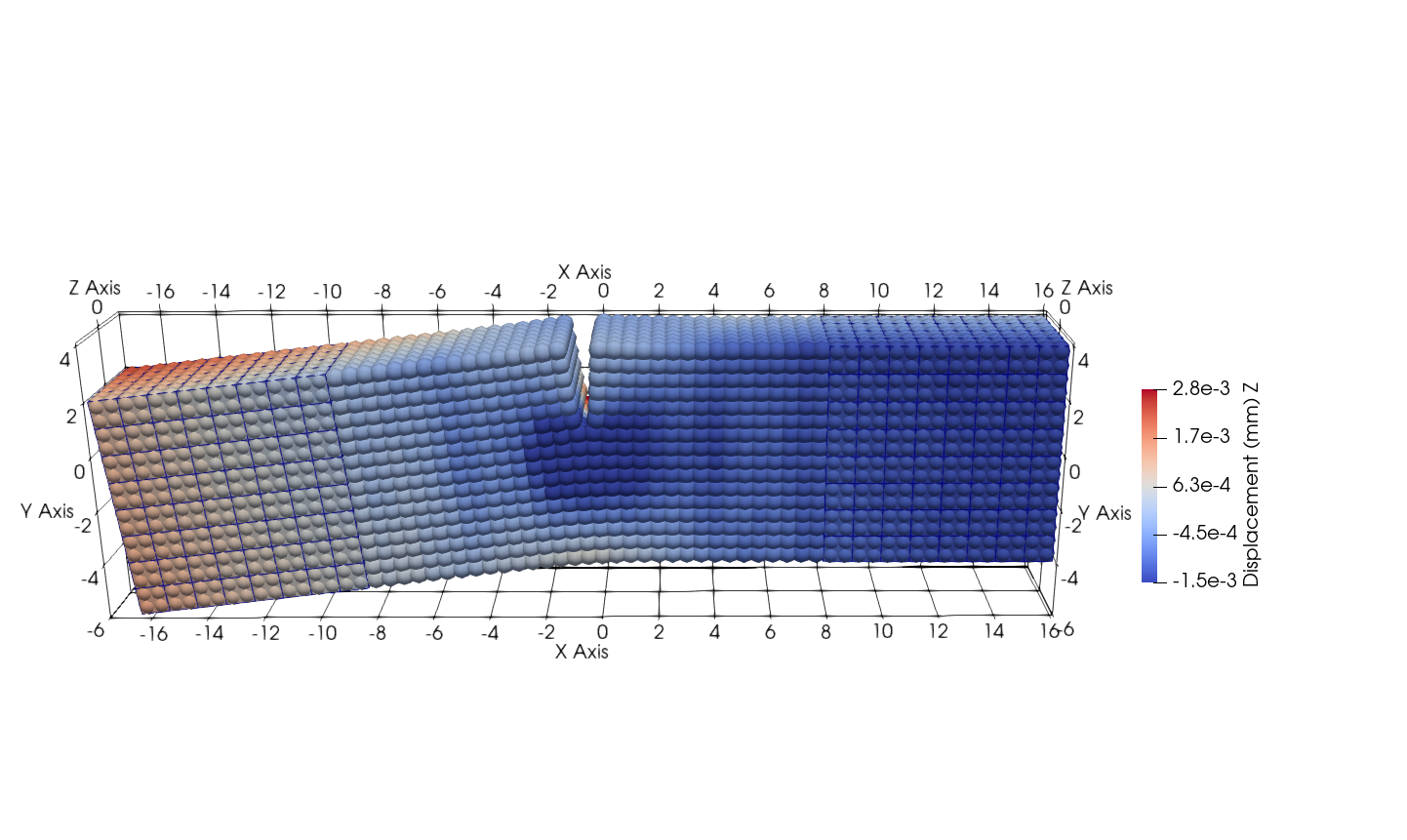}
     \subcaption{Displacement in $z$}\label{fig::DirichletNeumannBarZ}
     }
     \end{subfigure}

    \end{center}
    \caption{Bar with Neumann-Dirichlet Boundary Conditions. Displacements are magnified by a factor of $15$ to emphasize displacements.}\label{fig::NeumannBarImages}
\end{figure}

\begin{figure}
    \begin{center}
  \begin{subfigure}{0.32\textwidth}
    {   \centering
     \includegraphics[width=0.9\textwidth,trim=11.0cm 2.3cm 12.6cm 1.6cm,clip]{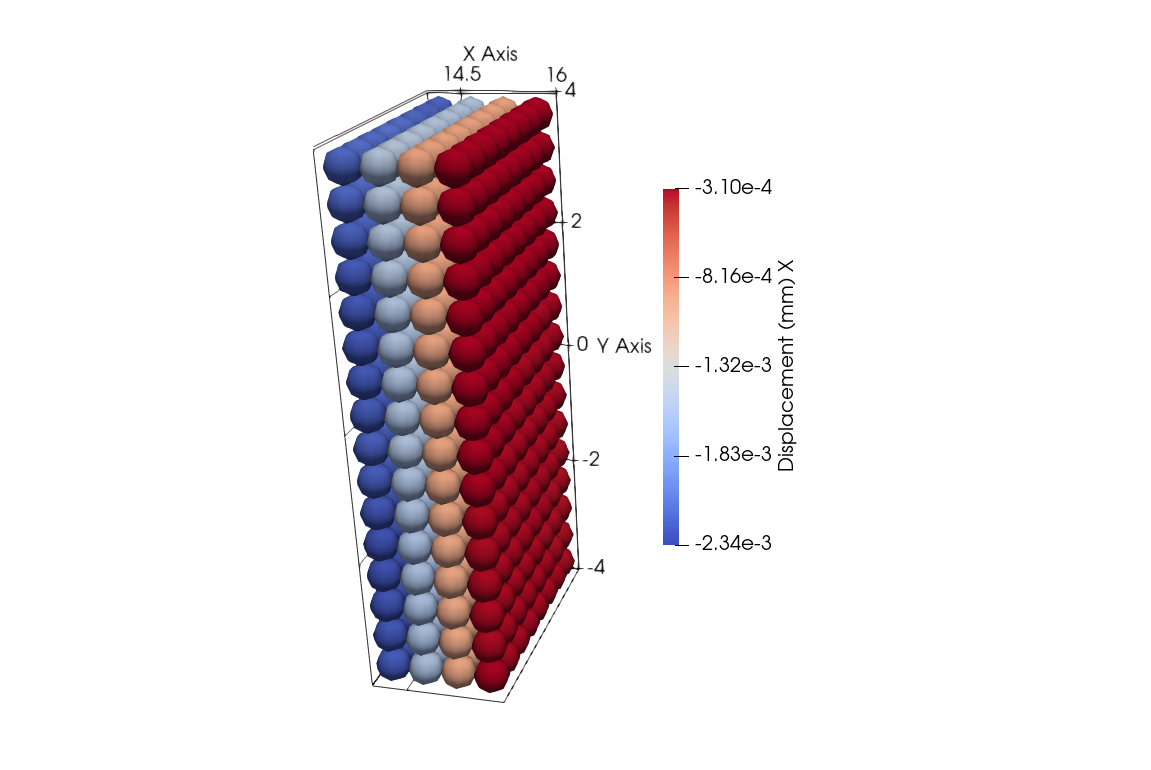}
     \subcaption{Displacement in $x$}\label{fig::DirichletNeumannBarDirichletBoundaryX}
     }
     \end{subfigure}
  \begin{subfigure}{0.32\textwidth}
    {   \centering
     \includegraphics[width=0.9\textwidth,trim=11.0cm 2.3cm 12.6cm 1.6cm,clip]{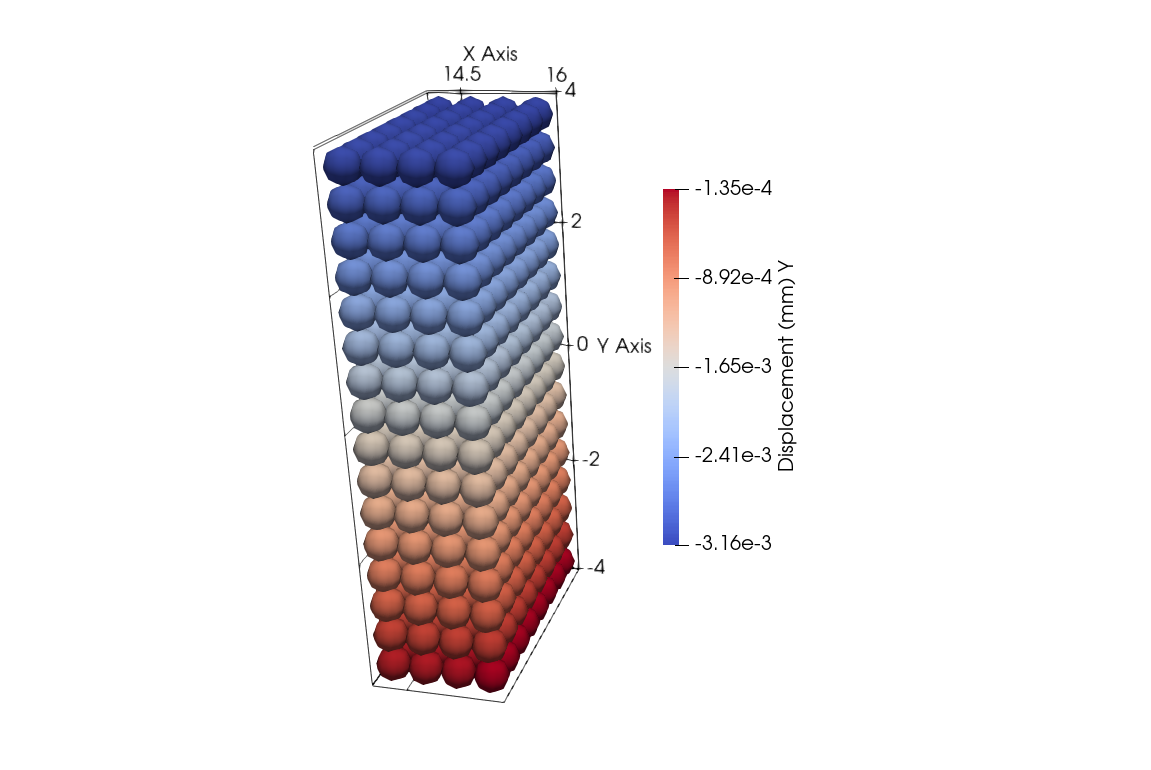}
     \subcaption{Displacement in $y$}\label{fig::DirichletNeumannBarDirichletBoundaryY}
     }
     \end{subfigure}
  \begin{subfigure}{0.32\textwidth}
    {   \centering
     \includegraphics[width=0.9\textwidth,trim=11.0cm 2.3cm 12.6cm 1.6cm,clip]{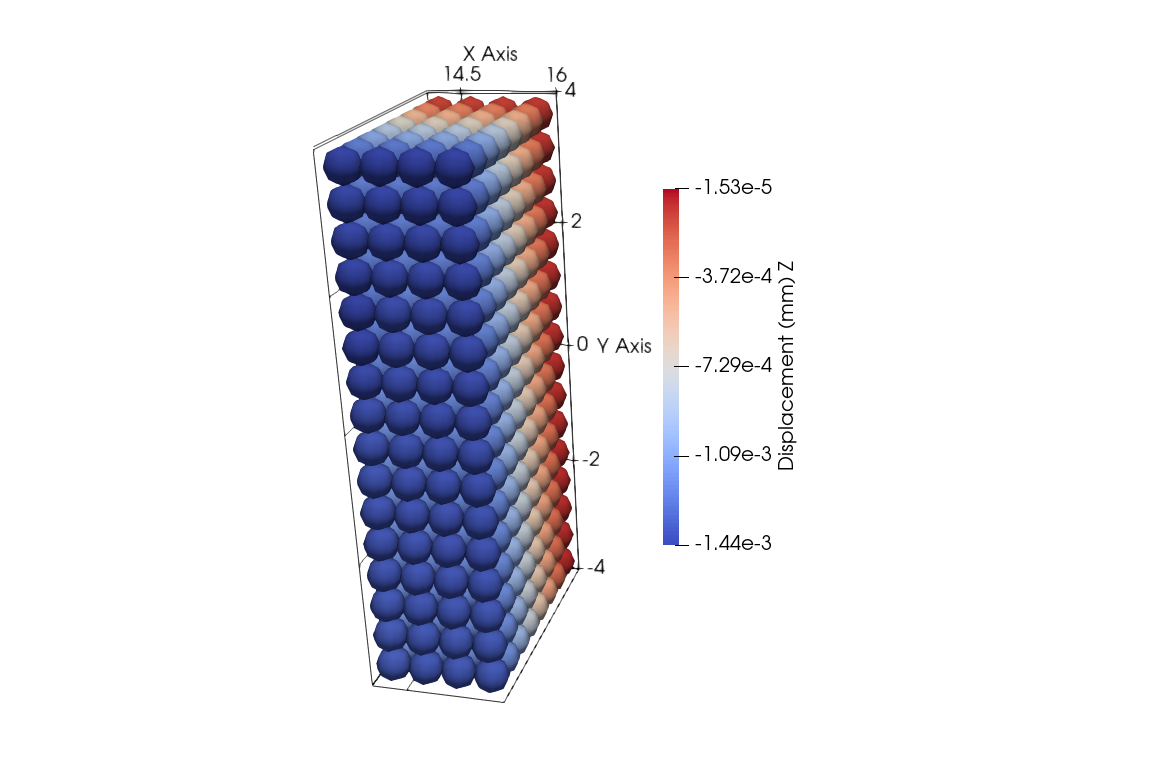}
     \subcaption{Displacement in $z$}\label{fig::DirichletNeumannBarDirichletBoundaryZ}
     }
     \end{subfigure}
     
   \begin{subfigure}{0.32\textwidth}
    {   \centering
     \includegraphics[width=0.9\textwidth,trim=6.5cm 2.9cm 10.4cm 2.4cm,clip]{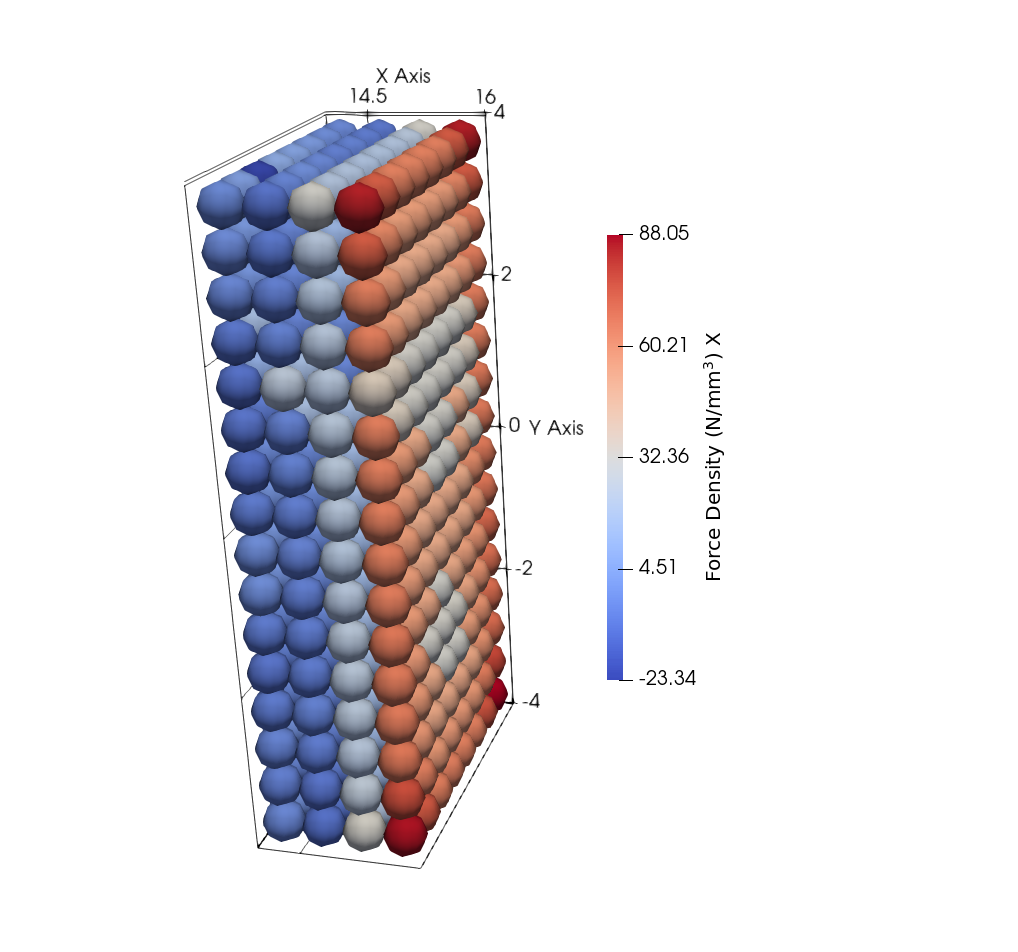}
     \subcaption{Force density in $x$}\label{fig::DirichletNeumannBarDirichletBoundaryXForce}
     }
     \end{subfigure}
  \begin{subfigure}{0.32\textwidth}
    {   \centering
     \includegraphics[width=0.9\textwidth,trim=6.5cm 2.9cm 10.4cm 2.4cm,clip]{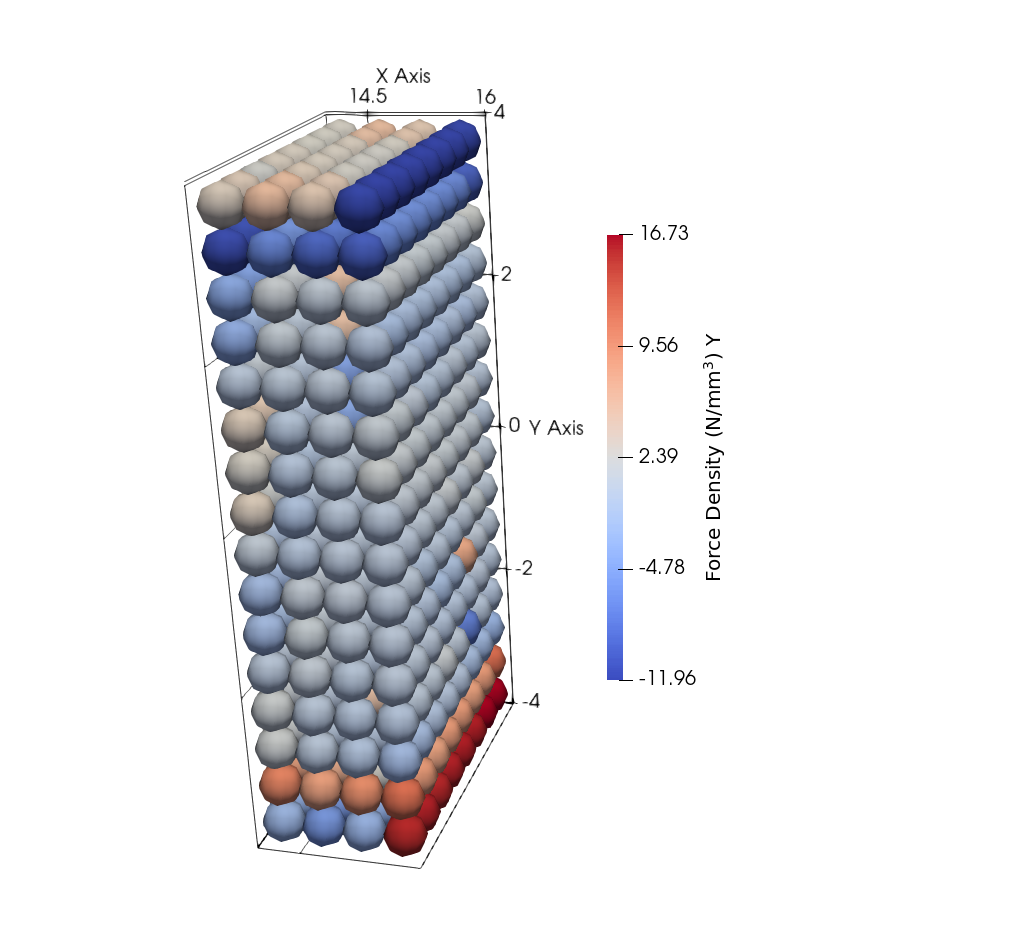}
     \subcaption{Force density in $y$}\label{fig::DirichletNeumannBarDirichletBoundaryYForce}
     }
     \end{subfigure}
  \begin{subfigure}{0.32\textwidth}
    {   \centering
     \includegraphics[width=0.9\textwidth,trim=6.5cm 2.9cm 10.4cm 2.4cm,clip]{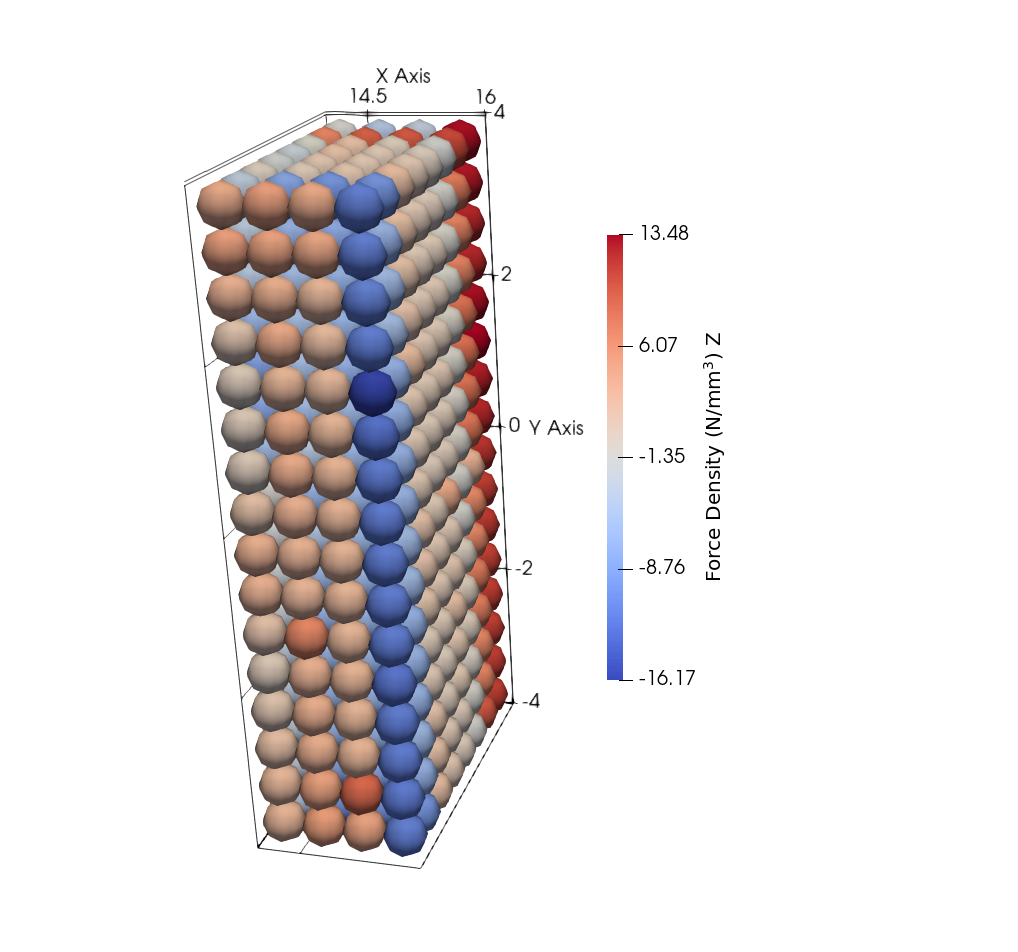}
     \subcaption{Force density in $z$}\label{fig::DirichletNeumannBarDirichletBoundaryZForce}
     }
     \end{subfigure}    
     
    \end{center}
    \caption{Nonlocal boundary conditions corresponding to the local Dirichlet boundary condition.}\label{fig::DirichletNeumannBarDirichletBoundary}
\end{figure}

\begin{figure}
    \begin{center}
  \begin{subfigure}{0.32\textwidth}
    {   \centering
     \includegraphics[width=0.9\textwidth,trim=10.4cm 3.5cm 12.6cm 1.4cm,clip]{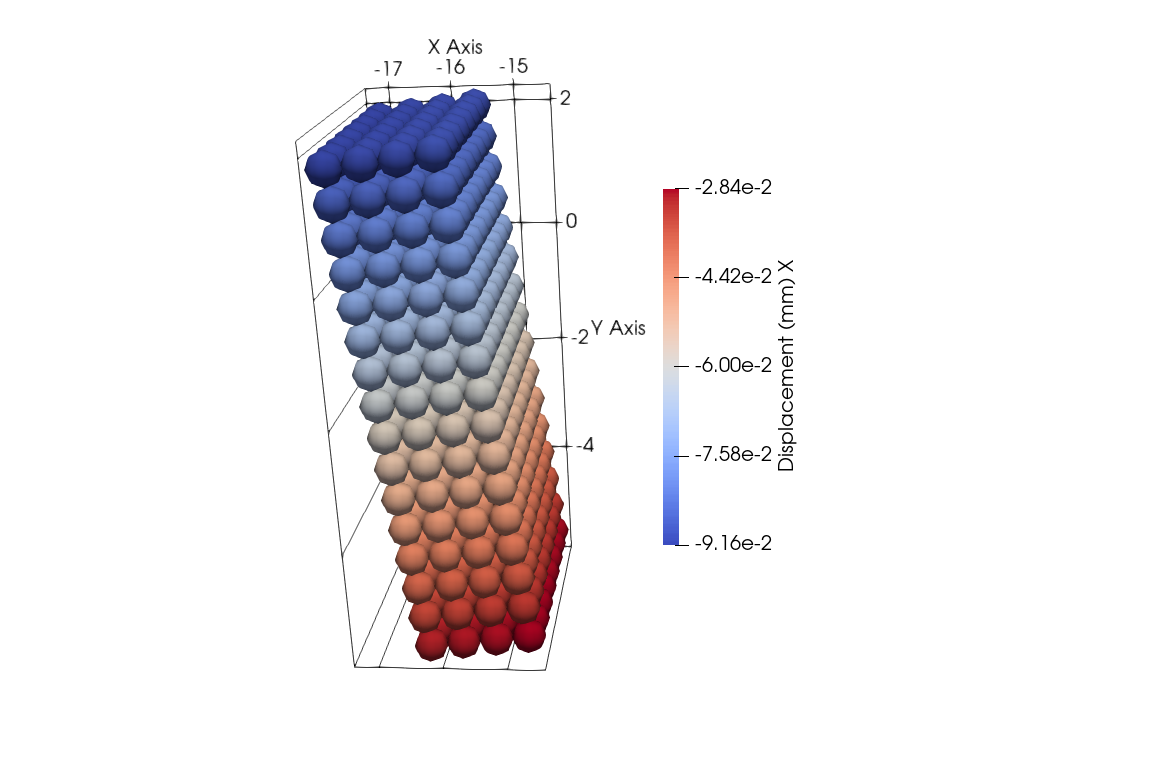}
     \subcaption{Displacement in $x$}\label{fig::DirichletNeumannBarNeumannBoundaryX}
     }
     \end{subfigure}
  \begin{subfigure}{0.32\textwidth}
    {   \centering
     \includegraphics[width=0.9\textwidth,trim=10.4cm 3.5cm 12.6cm 1.4cm,clip]{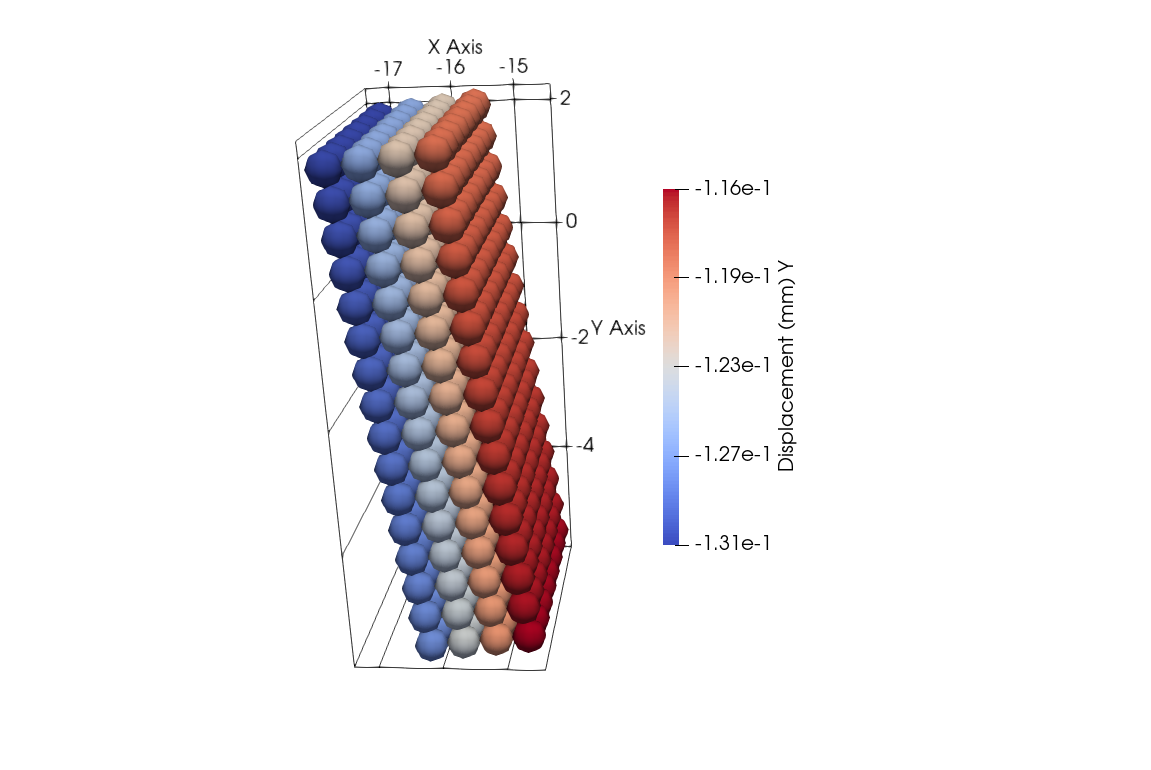}
     \subcaption{Displacement in $y$}\label{fig::DirichletNeumannBarNeumannBoundaryY}
     }
     \end{subfigure}
  \begin{subfigure}{0.32\textwidth}
    {   \centering
     \includegraphics[width=0.9\textwidth,trim=10.4cm 3.5cm 12.6cm 1.4cm,clip]{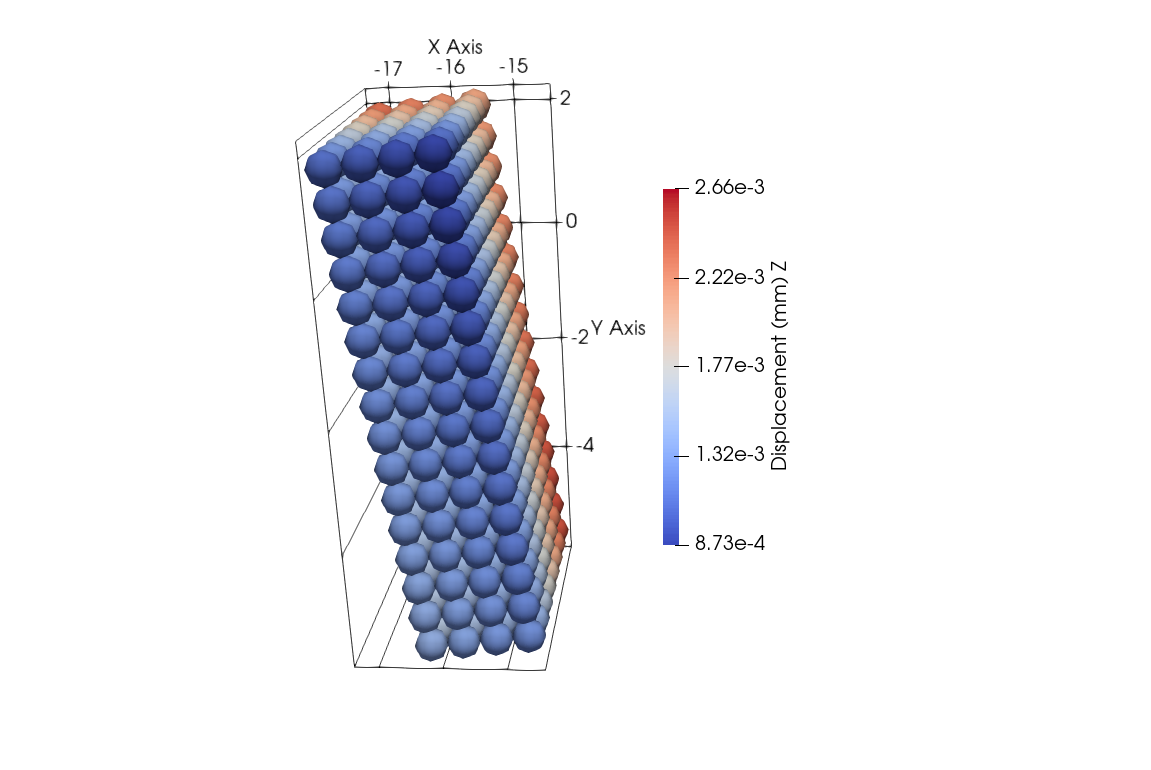}
     \subcaption{Displacement in $z$}\label{fig::DirichletNeumannBarNeumannBoundaryZ}
     }
     \end{subfigure}
     
  \begin{subfigure}{0.32\textwidth}
    {   \centering
     \includegraphics[width=0.9\textwidth,trim=5.7cm 4.3cm 11.0cm 2.1cm,clip]{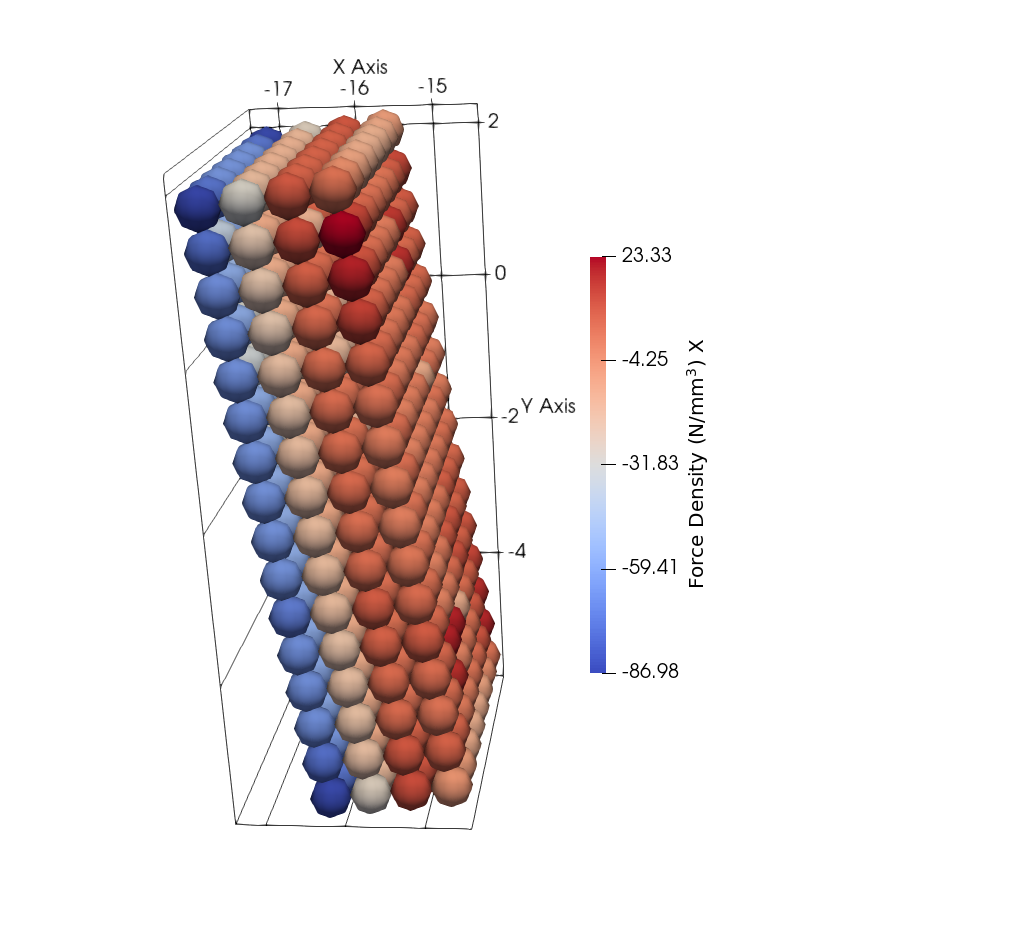}
     \subcaption{Force density in $x$}\label{fig::DirichletNeumannBarNeumannBoundaryXForce}
     }
     \end{subfigure}
  \begin{subfigure}{0.32\textwidth}
    {   \centering
     \includegraphics[width=0.9\textwidth,trim=5.7cm 4.3cm 11.0cm 2.1cm,clip]{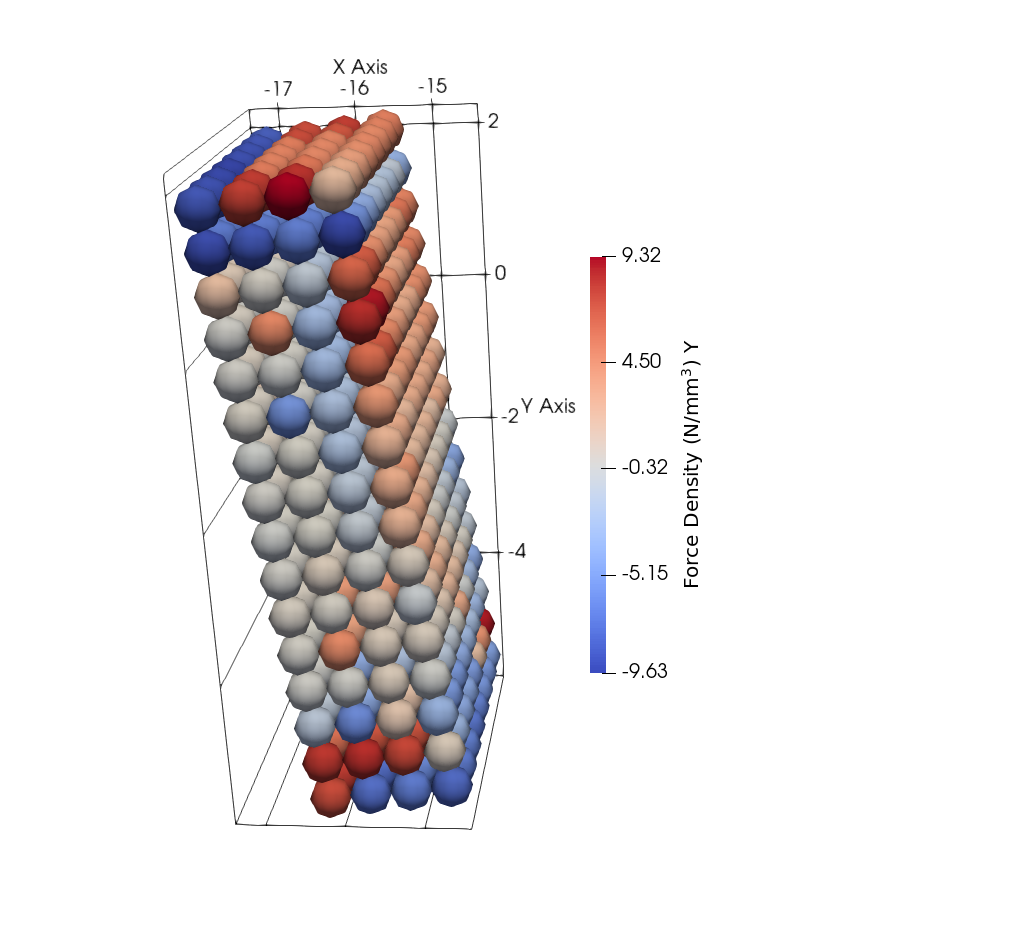}
     \subcaption{Force density in $y$}\label{fig::DirichletNeumannBarNeumannBoundaryYForce}
     }
     \end{subfigure}
  \begin{subfigure}{0.32\textwidth}
    {   \centering
     \includegraphics[width=0.9\textwidth,trim=5.7cm 4.3cm 11.0cm 2.1cm,clip]{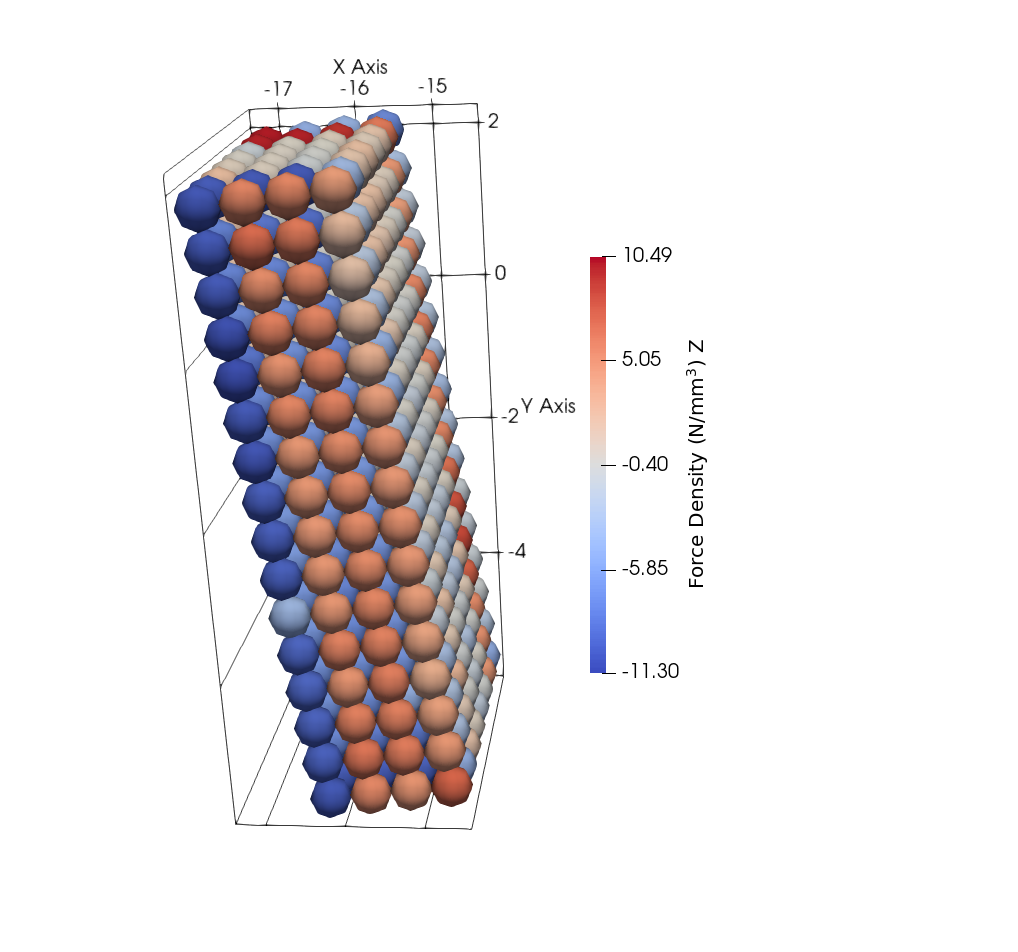}
     \subcaption{Force density in $z$}\label{fig::DirichletNeumannBarNeumannBoundaryZForce}
     }
     \end{subfigure}

    \end{center}
    \caption{Nonlocal boundary conditions corresponding to the local Neumann boundary condition.}\label{fig::DirichletNeumannBarNeumannBoundary}
\end{figure}

The results for this experiment are shown in Figures~\ref{fig::NeumannBarImages}, \ref{fig::DirichletNeumannBarDirichletBoundary}, and \ref{fig::DirichletNeumannBarNeumannBoundary}.  Displacement results for both the local and nonlocal dommains are shown in Figure~\ref{fig::NeumannBarImages} and are consistent with the expected behavior.  Displacement and force density results for the nonlocal model in the direct vicinity of the boundary conditions at the ends of the bar are presented in Figures~\ref{fig::DirichletNeumannBarDirichletBoundary} and \ref{fig::DirichletNeumannBarNeumannBoundary}.  The data in Figure~\ref{fig::DirichletNeumannBarDirichletBoundary} correspond to the end of the bar at $x=$ \SI{16.0}{\milli\meter}, adjacent to the local Dirichlet boundary conditions.  We see that the displacements in the nonlocal model vary approximately linearly over this region.  This is an intuitive result, and it suggests that a fair approximation of the prescribed displacement volume constraints for a comparable purely nonlocal model could be determined simply by assuming a linear displacement field in this region.  It is noted that the force density results in Figure~\ref{fig::DirichletNeumannBarDirichletBoundary} reflect the influence of the applied loading, the free surfaces, and the fixed displacement boundary conditions applied to the local model to prevent rigid body modes.

Results in Figure~\ref{fig::DirichletNeumannBarNeumannBoundary} show the displacements and force densities for the nonlocal model directly adjacent to the traction boundary condition that is applied to the local model.  The displacement field is again approximately linear, which is consistent with the expected behavior.  The force density field is much more complex, however, which is a reflection of the nonlocal interactions between material points in the peridynamic model.  These interactions are dictated by the nonlocal material model, the geometry (e.g., free surfaces), the discretization, and the value of the horizon.  A goal of the current study is to illustrate the complexity of determining prescribed force density volume constraints for a purely nonlocal model that approximate the effect of traction boundary conditions in a comparable local model, and to provide a systematic approach for doing so via LtN coupling.

\section{Summary}

In this study, we presented an optimization-based approach for local-to-nonlocal coupling and demonstrated its effectiveness on a series of computational examples.  The objective function for the optimization problem is the difference between local and nonlocal solutions in an overlap region, the controls are virtual boundary conditions on the local and nonlocal models at the coupling interface, and the models themselves act as constraints on the optimization problem.  Numerical examples on statics problems, carried out by coupling the \emph{Peridigm} and \emph{Albany} codes, demonstrate the accuracy of the coupling method, its convergence behavior, and its application to the geometry of a compact tension test experiment.  The coupling method is shown to be a viable means for carrying out simulations in which a nonlocal model is applied in the vicinity of a discontinuity, such as a crack, and a less computationally expensive local model is applied elsewhere.  The coupling approach is also shown to enable the conversion of local boundary conditions, including tractions, to prescribed displacements or force densities for the definition of nonlocal volume constraints.  Future work includes application of the coupling scheme for quasi-static problems that include stable crack growth, and an investigation of its use for approximating nonlocal volume constraints in transient dynamics simulations.

\section*{Acknowledgments}
M. D'Elia, P. Bochev and M. Perego were partially supported by the U.S. Department of Energy, Office of Advanced Scientific Computing Research under the Collaboratory on Mathematics and Physics-Informed Learning Machines for Multiscale and Multiphysics Problems (PhILMs) project. M. D'Elia and J. Trageser were also supported by the Sandia National Laboratories Laboratory-directed Research and Development (LDRD) program, project 218318.

Sandia National Laboratories is a multimission laboratory managed and operated by National Technology and Engineering Solutions of Sandia, LLC., a wholly owned subsidiary of Honeywell International, Inc., for the U.S. Department of Energy's National Nuclear Security Administration contract number DE-NA0003525. This paper, SAND2021-12618 R, describes objective technical results and analysis. Any subjective views or opinions that might be expressed in the paper do not necessarily represent the views of the U.S. Department of Energy or the United States Government.

\bibliographystyle{plain}
\bibliography{OBC}

\end{document}